\documentclass[12pt]{article}
\usepackage{latexsym, amssymb, amsmath, amscd, amsfonts, epsfig, graphicx, colordvi,verbatim,ifpdf,extarrows}
\usepackage{amsfonts, amsmath, amssymb}
\usepackage{amssymb,amsfonts,amsmath,latexsym,epsfig,cite, psfrag,eepic,color}
\usepackage{amscd,graphics}
\usepackage{latexsym, amssymb,  amsmath,amscd, amsfonts, epsfig, graphicx, colordvi,amsthm}

\usepackage{graphicx}
\usepackage{epstopdf}
\usepackage{color}
\usepackage{ifpdf}
\usepackage{fancybox}
\usepackage[font=small,labelfont=bf,labelsep=none]{caption}
\usepackage{float}

\allowdisplaybreaks

\usepackage[latin1]{inputenc}

\newtheorem{theorem}{Theorem}[section]
\newtheorem{prop}[theorem]{Proposition}

\theoremstyle{remark}

\theoremstyle{definition}
\newtheorem{definition}[theorem]{Definition}
\newtheorem{lemma}[theorem]{Lemma}
\newtheorem{core}[theorem]{Corollary}

\def\pf{\noindent{\it Proof.} }
\def\qed{\nopagebreak\hfill{\rule{4pt}{7pt}}
\medbreak}

\numberwithin{equation}{section}

\def\qed{\nopagebreak\hfill{\rule{4pt}{7pt}}
\medbreak}

\newlength{\boxedparwidth}
  {\begin{center} \begin{tabular}{|@{\hspace{.315in}}c@{\hspace{.15in}}|}
                  \hline \\ \begin{minipage}[t]{\boxedparwidth}
                  \setlength{\parindent}{.25in}}%
  {\end{minipage} \\ \\ \hline \end{tabular} \end{center}}

\begin{document}

\begin{center}
{\Large \bf Separable overpartition classes and excludant sizes of an overpartition}
\end{center}

\begin{center}
 {Y.H. Chen}$^{1}$, {Y.Q. Chen}$^{2}$, {Thomas Y. He}$^{3}$, {H.X. Huang}$^{4}$ and {X. Zhang}$^{5}$ \vskip 2mm

$^{1,2,3,4,5}$ School of Mathematical Sciences, Sichuan Normal University, Chengdu 610066, P.R. China

   \vskip 2mm

  $^1$chenyh@stu.sicnu.edu.cn, $^2$chenyuqian11@stu.sicnu.edu.cn, $^3$heyao@sicnu.edu.cn,  $^4$huanghaoxuan@stu.sicnu.edu.cn, $^5$zhangxi@stu.sicnu.edu.cn
\end{center}

\vskip 6mm   {\noindent \bf Abstract.}  An overpartition  is a partition such that the first occurrence (equivalently, the last occurrence) of a number may be overlined. In this article, we investigate three contents of overpartitions.
We first consider the $r$-chain minimal and maximal excludant sizes of an overpartition. Then, we study the second minimal excludant and mex sequence of an overpartition. Finally, we introduce $L_k$-overpartitions and $F_k$-overpartitions, which are separable overpartition classes.


\section{Introduction}

A partition $\pi$ of a positive integer $n$ is a finite non-increasing sequence of positive integers $\pi=(\pi_1,\pi_2,\ldots,\pi_\ell)$ such that $\pi_1+\pi_2+\cdots+\pi_\ell=n$. The empty sequence forms the only partition of zero. The $\pi_i$ are called the parts of $\pi$. Let $\ell(\pi)$ be the number of parts of $\pi$. The weight of $\pi$ is the sum of parts, denoted $|\pi|$. We use $\mathcal{P}(n)$ to denote the set of all partitions of $n$.

An overpartition, introduced by Corteel and Lovejoy \cite{Corteel-Lovejoy-2004},  is a partition such that the first occurrence (equivalently, the last occurrence) of a number may be overlined.
For example, there are fourteen overpartitions of $4$ in which the last occurrence of a number may be overlined.
\[(4),(\overline{4}),(3,1),(\overline{3},1),(3,\overline{1}),(\overline{3},\overline{1}),(2,2),(2,\overline{2}),
\]
\[(2,1,1),(\overline{2},1,1),(2,1,\overline{1}),(\overline{2},1,\overline{1}),(1,1,1,1),(1,1,1,\overline{1}).\]

We use $\overline{\mathcal{P}}$ to denote the set of all overpartitions and use $\overline{\mathcal{P}}(n)$ to denote the set of all overpartitions of $n$.
Let $\pi$ be an overpartition. We use $f_t(\pi)$ and $f_{\overline{t}}(\pi)$  to denote the number of parts equal to $t$ and $\overline{t}$ in $\pi$ respectively. We use $\ell_o(\pi)$ to denote the number of overlined parts in $\pi$. For a part $\pi_i$ of $\pi$, we say that $\pi_i$ is of size $t$ if $\pi_i=t$ or $\overline{t}$, denoted $|\pi_i|=t$. We adopt the following convention: For positive integer $c$ and nonnegative integer $d$, we define $\overline{c}+d=\overline{c+d}$.

The main objective of this article is to investigate two contents of overpartitions:
 \begin{itemize}
\item[(1)] the $r$-chain minimal and maximal excludant sizes of an overpartition;

\item[(2)] the second minimal excludant and mex sequence of an overpartition;

\item[(3)] $L_k$-overpartitions and $F_k$-overpartitions.

 \end{itemize}

 Throughout this article, we assume that $k$ and $r$ are positive integers.  We assume that $|q|<1$ and use the standard notation \cite{Andrews-1976}:
\[(a;q)_\infty=\prod_{i=0}^{\infty}(1-aq^i)\quad\text{and}\quad(a;q)_n=\frac{(a;q)_\infty}{(aq^n;q)_\infty}.\]

\subsection{Excludant sizes of an overpartition}

In  \cite{Andrews-Newman-2019}, Andrews and Newman undertook a combinatorial study of the minimal excludant of a partition, which was earlier introduced by Grabner and Knopfmacher \cite{Grabner-Knopfmacher-2006} under the name ``smallest gap". The minimal excludant of a partition $\pi$ is the smallest positive integer that is not a part of $\pi$, denoted  $mex(\pi)$. For $n\geq 0$, they defined
 \[
\sigma mex(n)=\sum_{\pi \in \mathcal{P}(n)}{mex}(\pi),
\]
and proved the following theorem.
\begin{theorem}\cite[Theorem 1]{Andrews-Newman-2019} For $n\geq 0$,
\[\sigma mex(n)=D_2(n),\]
where $D_2(n)$ is the number of partitions of $n$ into distinct parts using two colors.
\end{theorem}

Several generalizations of this concept have been investigated. In \cite{Bhoria-Eyyunni-Maji-2023}, Bhoria, Eyyunni and Maji  generalized the minimal excludant of a partition and studied the $r$-chain minimal excludant of a partition. The $r$-chain minimal excludant of a partition $\pi$ is the smallest positive integer $t$ such that none of $t,t+1,\ldots,t+r-1$ appears as a part of $\pi$, denoted $mex(\pi;r)$. They obtained the following result.
\begin{theorem}\cite[Theorem 2.6]{Bhoria-Eyyunni-Maji-2023}
\begin{equation*}
\sum_{n=0}^\infty\sum_{\pi\in\mathcal{P}(n)}mex(\pi;r)q^n=-\frac{r-1}{(q;q)_\infty}+\frac{(q^{r+1};q^{r+1})_\infty}{(q;q)_\infty}\sum_{m=1}^r\frac{1}{(q^m;q^{r+1})_\infty}.
\end{equation*}
\end{theorem}

Inspired by the minimal excludant of a partition, Chern \cite{Chern-2021} investigated the maximal excludant of a partition. The maximal excludant of a partition $\pi$ is the largest nonnegative integer smaller than the largest part of $\pi$ that does not appear in $\pi$. In \cite{Bhoria-Eyyunni-Li-2024}, Bhoria, Eyyunni and Li introduced the $r$-chain maximal excludant of a partition, which is a generalization of the maximal excludant of a partition analogous to the $r$-chain minimal excludant of a partition. The $r$-chain maximal excludant $maex(\pi;r)$ of a partition $\pi$ is the largest positive integer $t$ less than the largest part of $\pi$ such that the positive integers $t,t-1,\ldots,t-r+1$ do not occur as parts in $\pi$. If there is no such positive integer $t$ for $\pi$, then we set $maex(\pi;r)=0$. They obtained the following result.
\begin{theorem}\cite[Theorem 1.10]{Bhoria-Eyyunni-Li-2024}
\begin{equation*}
\sum_{n=0}^\infty\sum_{\pi\in\mathcal{P}(n)}\left(\ell(\pi)-maex(\pi;r)+\Omega_r(\pi)\right)q^n=\frac{(q^{r+1};q^{r+1})_\infty}{(q;q)_\infty}+\frac{1}{(q;q)_\infty}\sum_{n=1}^\infty\frac{q^n(q^{r+1};q^{r+1})_n}{1-q^n}.
\end{equation*}
\end{theorem}
Here, $\Omega_r(\pi)$ is defined in \cite[(1.3)]{Bhoria-Eyyunni-Li-2024}.

In \cite[Definition 1]{Kaur-Rana-Eyyunni-2024}, Kaur, Rana and Eyyunni considerer the second minimal excludant of a partition, which is the second smallest integer  missing from an integer partition $\pi$, denoted  $mex_{2}(\pi)$. For $n\geq 0$, they defined
\[
\sigma_2{mex}(n)= \sum_{\pi \in \mathcal{P}(n)}{mex}_2(\pi),
\]
and obtained the following theorem.

\begin{theorem}\cite[Theorem 4]{Kaur-Rana-Eyyunni-2024}
The generating function of $\sigma_2{mex}(n)$ is
\[
\sum_{n=0}^\infty\sigma_2mex(n)q^n=\frac{1}{(q;q)_\infty}\left(\frac{1}{1-q}-\sum_{s=0}^\infty(s-1)q^{\binom{s+1}{2}}\right).
\]
\end{theorem}

For $n\geq 0$ and $t\geq 1$, let ${\Delta}_t(n)$ be the number of overpartitions $\pi$ of $n$ satisfying ${\text{mex}}_2(\pi)-{\text{mex}}(\pi)=t$. Kaur, Rana and Eyyunni gave the generating function of  ${\Delta}_t(n)$.

\begin{theorem}\cite[Theorem 6]{Kaur-Rana-Eyyunni-2024}
For $t\geq1$,
\[
\sum_{n=0}^\infty{\Delta}_t(n)q^n=\frac{q^{t-1}}{(q^2;q)_\infty}\sum_{r=0}^\infty q^{\binom{r+1}{2}}.
\]
\end{theorem}

Kaur, Rana and Eyyunni \cite[Definition 2]{Kaur-Rana-Eyyunni-2024} also introduced the mex sequence of a partition, which is defined as the longest sequence of consecutive missing integers starting from the minimal excludant of a partition. For $n\geq 0$, let ${p}_r^{mex}(n)$ denote the number of partitions of $n$ whose mex sequence has length at least $r$.
Kaur, Rana and Eyyunni showed the following result.

\begin{prop}\cite[Proposition 9]{Kaur-Rana-Eyyunni-2024}\label{prop-add-1}
For $n\geq 0$, ${p}_r^{mex}(n)=q(n)$ if and only if $r\geq n$, where $q(n)$ is the number of distinct partitions of $n$.
\end{prop}

The generating function of ${p}_r^{mex}(n)$ was given by Kaur, Rana and Eyyunni.

\begin{theorem}\cite[Theorem 10]{Kaur-Rana-Eyyunni-2024}
\[
\sum_{n=0}^\infty{p}_r^{mex}(n)q^n=\frac{1}{(q;q^2)_\infty(q^{r+1};q^2)_\infty}.
\]
\end{theorem}

In \cite{Dhar-Mukhopadhyay-Sarma-2023,Dhar-Mukhopadhyay-Sarma-2024}, Dhar, Mukhopadhyay and Sarma extended the minimal excludant of partitions to overpartitions.
They defined
the minimal excludant $\overline{mex}(\pi)$ of an overpartition $\pi$ as the smallest positive integer such that there are no parts of size $\overline{mex}(\pi)$ in  $\pi$. The $mes_{1,1,1}(\pi)$ in \cite{He-Huang-Li-Zhang-2024} is $\overline{mex}(\pi)$.
In the remaining of this subsection, an overpartition is defined as a partition in which the last occurrence of a number may be overlined.
For example, the minimal excludants of overpartitions of $4$ are given in the following tables.
\begin{center}
\begin{tabular}{|c|c|c|c|c|c|c|c|c|}
  \hline
  $\pi$&$(4)$&$(\overline{4})$&$(3,1)$&$(\overline{3},1)$&$(3,\overline{1})$&$(\overline{3},\overline{1})$&$(2,2)$&$({2},\overline{2})$\\
  \hline
  $\overline{mex}(\pi)$&$1$&$1$&$2$&$2$&$2$&$2$&$1$&$1$\\
  \hline
\end{tabular}
\end{center}

\begin{center}
\begin{tabular}{|c|c|c|c|c|c|c|}
  \hline
  $\pi$&$(2,1,1)$&$(\overline{2},1,1)$&$(2,{1},\overline{1})$&$(\overline{2},{1},\overline{1})$&$(1,1,1,1)$&$({1},1,1,\overline{1})$\\
  \hline
  $\overline{mex}(\pi)$&$3$&$3$&$3$&$3$&$2$&$2$\\
  \hline
\end{tabular}
\end{center}

In this article, we extend the $r$-chain minimal excludant, the $r$-chain maximal excludant, the second minimal excludant and mex sequence of a partition to an overpartition.

{\noindent \bf The $r$-chain minimal and maximal excludant sizes of an overpartition}

\begin{definition}[$r$-chain minimal excludant size]\label{defi-min}
The $r$-chain minimal excludant size of an overpartition $\pi$, denoted $mes(\pi;r)$, is the smallest positive integer $t$ such that there are no parts of size $t,t+1,\ldots,t+r-1$ in $\pi$.
\end{definition}

For example, the $2$-chain minimal excludant sizes of  overpartitions of $4$ are given in the following tables.

\begin{center}
\begin{tabular}{|c|c|c|c|c|c|c|c|c|}
  \hline
  $\pi$&$(4)$&$(\overline{4})$&$(3,1)$&$(\overline{3},1)$&$(3,\overline{1})$&$(\overline{3},\overline{1})$&$(2,2)$&$(2,\overline{2})$\\
  \hline
  ${mes}(\pi,2)$&$1$&$1$&$4$&$4$&$4$&$4$&$3$&$3$\\
  \hline
\end{tabular}
\end{center}

\begin{center}
\begin{tabular}{|c|c|c|c|c|c|c|}
  \hline
  $\pi$&$(2,1,1)$&$(\overline{2},1,1)$&$(2,1,\overline{1})$&$(\overline{2},1,\overline{1})$&$(1,1,1,1)$&$(1,1,1,\overline{1})$\\
  \hline
  ${mes}(\pi,2)$&$3$&$3$&$3$&$3$&$2$&$2$\\
  \hline
\end{tabular}
\end{center}

\begin{definition}[$r$-chain maximal excludant size]\label{defi-max}
The $r$-chain maximal excludant size of an overpartition $\pi$, denoted $maes(\pi;r)$, is the largest positive integer $t$ less than the size of the largest part of $\pi$ such that $t\geq r$ and there are no parts of size $t,t-1,\ldots,t-r+1$  in $\pi$. If there is no such positive integer $t$ for $\pi$, then we set $maes(\pi;r)=0$.
\end{definition}

For example, the overpartitions $\pi$ of $6$ with positive $2$-chain maximal excludant sizes are given in the following tables.

\begin{center}
\begin{tabular}{|c|c|c|c|c|c|c|}
  \hline
  $\pi$&$(6)$&$(\overline{6})$&$(5,1)$&$(\overline{5},1)$&$(5,\overline{1})$&$(\overline{5},\overline{1})$\\
  \hline
  ${maes}(\pi,2)$&$5$&$5$&$4$&$4$&$4$&$4$\\
  \hline
\end{tabular}
\end{center}

\begin{center}
\begin{tabular}{|c|c|c|c|c|c|c|}
  \hline
  $\pi$&$(4,1,1)$&$(\overline{4},1,1)$&$(4,1,\overline{1})$&$(\overline{4},1,\overline{1})$&$(3,3)$&$(3,\overline{3})$\\
  \hline
  ${maes}(\pi,2)$&$3$&$3$&$3$&$3$&$2$&$2$\\
  \hline
\end{tabular}
\end{center}

For $n\geq 0$, we define
\[
\sigma_r{mes}(n)=\sum_{\pi\in\overline{\mathcal{P}}(n)}{mes}(\pi,r),
\]
and
\[
\sigma_r{maes}(n)=\sum_{\pi\in\overline{\mathcal{P}}(n)}{maes}(\pi,r).
\]

We will give the generating functions of $\sigma_r{mes}(n)$ and $\sigma_r{maes}(n)$. For easier expression, we use the following notations. For a positive integer $t$, we set
\[\omega(t)=1+\sum_{i=1}^r2q^{it}.\]

For positive integers $m$ and $n$, we set
\[\omega(m;0)=1,\]
\[\omega(m;\infty)=\prod_{i=0}^\infty\omega(m+i),\]
and
\[\omega(m;n)=\prod_{i=0}^{n-1}\omega(m+i).\]

\begin{theorem}\label{add-thm-gen-m-a-es}
\begin{equation}\label{gen-mes-r}
\begin{split}
\sum_{n=0}^\infty\sigma_r{mes}(n)q^n&=\frac{(-q;q)_{\infty}}{(q;q)_{\infty}}+\sum_{n=1}^{\infty}\frac{2(q^n+2q^{2n}+\cdots+rq^{rn})}{\omega(n)}\\
&\qquad\qquad\qquad\times\left(\frac{(-q;q)_{\infty}}{(q;q)_{\infty}}-\omega(1;n-1)\frac{(-q^n;q)_{\infty}}{(q^n;q)_{\infty}}+\omega(1;\infty)\right),
\end{split}
\end{equation}
and
\begin{equation}\label{gen-maes-r}
\begin{split}
\sum_{n=0}^\infty\sigma_r{maes}(n)q^n&=r\left(\frac{(-q;q)_{\infty}}{(q;q)_{\infty}}-\omega(1;\infty)\right)+\frac{(-q;q)_{\infty}}{(q;q)_{\infty}}\sum_{n=1}^\infty\frac{2q^n}{1-q^{2n}}\\
&\qquad-\sum_{n=1}^\infty q^n\omega(1;n-1)\frac{(-q^{n+1};q)_{\infty}}{(q^n;q)_{\infty}}(1+w(n)).
\end{split}
\end{equation}
\end{theorem}

In \cite{He-Huang-Li-Zhang-2024}, He, Huang, Li and Zhang obtained the following generating function of $\sigma_1{mes}(n)$.
\begin{equation*}\label{gen-sigma-mes-eqn}
\sum_{n=0}^\infty\sigma_1{mes}(n)q^n
=\frac{(-q;q)_{\infty}}{(q;q)_{\infty}}\sum_{k=0}^{\infty}\frac{2^kq^{\binom{k+1}{2}}}{(-q;q)_k}.
\end{equation*}
Combining with  \eqref{gen-mes-r}, we can get
\begin{equation*}
\frac{(-q;q)_{\infty}}{(q;q)_{\infty}}\sum_{k=1}^{\infty}\frac{2^kq^{\binom{k+1}{2}}}{(-q;q)_k}=\sum_{n=1}^{\infty}\frac{2q^n}{1+2q^n}\left(\frac{(-q;q)_{\infty}}{(q;q)_{\infty}}-(-2q;q)_{n-1}\frac{(-q^n;q)_{\infty}}{(q^n;q)_{\infty}}+(-2q;q)_{\infty}\right).
\end{equation*}

{\bf \noindent The second minimal excludant and mex sequence of an overpartion}

\begin{definition}[The second minimal excludant]\label{Second Minimal Excludant}
The second minimal excludant $\overline{mex}_2(\pi)$ of an overpartition $\pi$  is the second smallest positive integer such that there are no parts of size $\overline{mex}_2(\pi)$ in $\pi$.
 \end{definition}

 For example, the second minimal excludants of  overpartitions of $4$ are given in the following tables.

\begin{center}
\begin{tabular}{|c|c|c|c|c|c|c|c|c|}
  \hline
  $\pi$&$(4)$&$(\overline{4})$&$(3,1)$&$(\overline{3},1)$&$(3,\overline{1})$&$(\overline{3},\overline{1})$&$(2,2)$&$({2},\overline{2})$\\
  \hline
  $\overline{mex}_2(\pi)$&$2$&$2$&$4$&$4$&$4$&$4$&$3$&$3$\\
  \hline
\end{tabular}
\end{center}

\begin{center}
\begin{tabular}{|c|c|c|c|c|c|c|}
  \hline
  $\pi$&$(2,1,1)$&$(\overline{2},1,1)$&$(2,{1},\overline{1})$&$(\overline{2},1,\overline{1})$&$(1,1,1,1)$&$(1,1,1,\overline{1})$\\
  \hline
  $\overline{mex}_2(\pi)$&$4$&$4$&$4$&$4$&$3$&$3$\\
  \hline
\end{tabular}
\end{center}

For $n\geq 0$, we define
\begin{equation*}\label{mex-1}
{\sigma}_2\overline{mex}(n)=\sum_{\pi \in \overline{\mathcal{P}}(n)}\overline{\text{mex}}_2(\pi).
\end{equation*}
For example, we have ${\sigma}_2\overline{mex}(4)=48$. We will prove the following result.
\begin{theorem}\label{mex}
\begin{equation*}
\sum_{n=0}^{\infty} {\sigma}_2\overline{mex}(n)q^n=\frac{(-q;q)_\infty}{(q;q)_\infty}\left(\frac{1}{1-q}+\sum_{s=0}^{\infty}\frac{(2-s)q^{\binom{s+1}{2}}2^{s-1}}{(-q;q)_s}\right).
\end{equation*}
\end{theorem}

For $n\geq 0$ and $t\geq1$ we define $\overline{\Delta}_t(n)$ to be the number of overpartitions $\pi$ of $n$ satisfying $\overline{mex}_2(\pi)-\overline{mex}(\pi)=t$. For instance, we have $\overline{\Delta}_1(4)=8$ and $\overline{\Delta}_2(4)=6$. We will show the following theorem.
\begin{theorem}\label{delta_1}
For $t\geq 1$, 
\begin{equation}\label{delta-2}
\sum_{n=0}^{\infty}\overline{\Delta}_t(n)q^n=(q^{t-1}-q^t)\frac{(-q;q)_\infty}{(q;q)_\infty}\sum_{m=0}^\infty\frac{q^{\binom{m+t}{2}}2^{m+t-1}}{(-q;q)_{m+t}}.
\end{equation}
\end{theorem}

Then, we turn to considering the mex sequence of an overpatition.
\begin{definition}[mex sequence]
The mex sequence of an overpatition is the longest sequence of consecutive missing sizes in the overpartition, starting from the size of its minimal excludant.
\end{definition}

For example, the following table shows the mex sequences of some overpartitions of $4$.
\begin{center}
\begin{tabular}{|c|c|c|c|c|c|c|c|c|}
  \hline
  $\text{overpartiton}$&$(4)$&$(\overline{4})$&$(3,1)$&$(\overline{3},1)$&$(3,\overline{1})$&$(\overline{3},\overline{1})$&$(2,2)$&$({2},\overline{2})$\\
  \hline
  $\text{mex sequence}$&$(1,2,3)$&$(1,2,3)$&$(2)$&$(2)$&$(2)$&$(2)$&$(1)$&$(1)$\\
  \hline
  $\text{length of mex sequence}$&$3$&$3$&$1$&$1$&$1$&$1$&$1$&$1$\\
  \hline
\end{tabular}
\end{center}

It is worth mentioning that the mex sequence of an overpartition can be infinitely long. For example, the  overpartitions
\[(2,1,1),(\overline{2},1,1),(2,{1},\overline{1}),(\overline{2},{1},\overline{1}),(1,1,1,1),({1},1,1,\overline{1}).\]
have such mex sequence.

For $n\geq 0$, let $\overline{p}_r^{mex}(n)$ denote the set of overpartitions of $n$ whose mex sequence has length at
 least $r$. With a similar argument in the proof of Proposition \ref{prop-add-1} in \cite{Kaur-Rana-Eyyunni-2024}, we can get the following result.
\begin{prop}
For $n\geq 0$, $\overline{p}_r^{mex}(n)=\overline{q}(n)$ if and only if $r\geq n$, where $\overline{q}(n)$ is the number of overpartitions $\pi$ of $n$ such that the sizes of parts in $\pi$ are distinct.
\end{prop}

We will give the generating function of $\overline{p}_r^{mex}(n)$.
\begin{theorem}\label{P_mex_1}
\begin{equation*}
\sum_{n=0}^\infty\overline{p}_r^{mex}(n)q^n=\sum_{m=0}^\infty \frac{(-1;q)_m}{(q;q)_m}q^{(r+1)m}(-2q^{m+1};q)_\infty.
\end{equation*}
\end{theorem}

For $n\geq 0$, let $\overline{p}_r(n)$ be the number of overpartitions $\pi$ of $n$ such that the number of parts of size $t$ in $\pi$ does not exceed $1$ for $t>m$, where $m$ is the largest integer such that the number of parts of size $m$ in $\pi$ is greater than $r$ with the convention that the number of parts of size $0$ in $\pi$ is $+\infty$.
For example, we have $\overline{p}_3(4)=8$ since there are the following eight overpartitions counted by $\overline{p}_3(4)$.
\[(4),(\overline{4}),(3,1),(\overline{3},1),(3,\overline{1}),(\overline{3},\overline{1}),(1,1,1,1),({1},1,1,\overline{1}).\]

It is clear that 
\[\sum_{n=0}^\infty\overline{p}_r(n)q^n=\sum_{m=0}^\infty \frac{(-1;q)_m}{(q;q)_m}q^{(r+1)m}(-2q^{m+1};q)_\infty.\]
As a corollary of Theorem \ref{P_mex_1}, we have the following result.
\begin{core}For $n\geq 0$,
\[\overline{p}_r^{mex}(n)=\overline{p}_r(n).\]
\end{core}

\subsection{$L_k$-overpartitions and $F_k$-overpartitions}

Andrews \cite{Andrews-2022} first introduced separable integer partition classes with modulus $k$ and analyzed some well-known theorems, such as the first G\"ollnitz-Gordon identity, Schur's partition theorem, partitions  with $n$ copies of $n$, and so on.

Based on separable integer partition classes with modulus $2$. Passary \cite[Section \ref{proof_3}]{Passary-2019} studied two cases of partitions with parts separated by parity, the little G\"ollnitz identities and the second G\"ollnitz-Gordon identity.

 In \cite{Chen-He-Tang-Wei-2024}, Chen, He, Tang and Wei  investigated the remaining six cases of partitions with parts separated by parity with the aid of separable integer partition classes with modulus $2$. They extended separable integer partition classes with modulus $1$ to overpartitions, called separable overpartition classes, and then studied overpartitions and the overpartition analogue of Rogers-Ramanujan identities from the view of separable overpartition classes.

By virtue of separable integer partition classes with modulus $k$, He, Huang, Li and Zhang \cite{He-Huang-Li-Zhang-2025} considered the partitions whose parts are congruent to $a$ or $b$ modulo $k$. They extended separable integer partition classes with modulus $k$ to overpartitions and then give the generating function for $(k,r)$-modulo overpartitions, which are the overpartitions satisfying certain congruence conditions.

We recall the definition of separable integer partition classes given by  Chen, He, Tang and Wei \cite{Chen-He-Tang-Wei-2024}, which is stated as follows.

 \begin{definition}[separable overpartition class]\label{defi-separable}
A separable overpartition class $\mathcal{P}$ is a subset of all the overpartitions satisfying the following{\rm:}

There is a subset $\mathcal{B}\subset\mathcal{P}$ {\rm(}$\mathcal{B}$ is called the basis of $\mathcal{P}${\rm)} such that for each integer $m\geq 1$, the number of overpartitions in $\mathcal{B}$ with $m$ parts is finite and every overpartition in $\mathcal{P}$ with $m$ parts is uniquely of the form
\begin{equation}\label{over-form-1}
(\lambda_1+\mu_1)+(\lambda_2+\mu_2)+\cdots+(\lambda_m+\mu_m),
\end{equation}
where $(\lambda_1,\lambda_2,\ldots,\lambda_m)$ is an overpartition in $\mathcal{B}$ and $(\mu_1,\mu_2,\ldots,\mu_m)$ is a  non-increasing sequence of nonnegative integers. Moreover, all overpartitions of the form \eqref{over-form-1} are in $\mathcal{P}$.
\end{definition}

Assume that $\mathcal{P}$ is a  separable overpartition class and $\mathcal{B}$ is the basis of $\mathcal{P}$. For $m,j\geq 1$, let $\mathcal{B}(m,j)$ (resp. $\mathcal{B}(m,\overline{j}))$ be the set of overpartitions in $\mathcal{B}$ with $m$ parts and the largest part being $j$ (resp. $\overline{j}$). Then, the generating function for  the overpartitions in $\mathcal{P}$ is
 \begin{align*}\label{useful}
  &\quad\sum_{\pi\in\mathcal{P}}z^{\ell_o(\pi)}q^{|\pi|}\\
  &=1+\sum_{s=1}^k\sum_{m= 1}^\infty\sum_{j=1}^\infty\frac{1}{(q;q)_{k(m-1)+s}}\left(\sum_{\lambda\in\mathcal{B}(k(m-1)+s,j)}z^{\ell_o(\lambda)}q^{|\lambda|}+\sum_{\lambda\in\mathcal{B}(k(m-1)+s,\overline{j})}z^{\ell_o(\lambda)}q^{|\lambda|}\right).
  \end{align*}

In this article, we will introduce two types of separable overpartition classes.
\begin{definition}\label{defi-L-k}
 Let $\pi=(\pi_1,\pi_2,\ldots,\pi_\ell)$ be an overpartition such that the last occurrence of a number may be overlined, and if $\pi_i$ is overlined then $\ell-i\equiv0\pmod{k}$. Then, we say that $\pi$ is a $L_k$-overpartition.
\end{definition}

For $k=1$, a $L_1$-overpartition is an overpartition in which the last occurrence  of a number may be overlined. For $k=2$, there are twelve $L_2$-overpartitions of $4$.
\[(4),(\overline{4}),(3,1),(3,\overline{1}),(2,2),(2,\overline{2}),(2,1,1),
\]
\[(\overline{2},1,1),(2,1,\overline{1}),(\overline{2},1,\overline{1}),(1,1,1,1),(1,1,1,\overline{1}).\]

\begin{definition}\label{defi-F-k}
Let $\pi=(\pi_1,\pi_2,\ldots,\pi_\ell)$ be an overpartition such that the first occurrence of a number may be overlined, and if $\pi_i$ is overlined then $\ell-i\equiv-1\pmod{k}$. Then, we say that $\pi$ is a $F_k$-overpartition.
\end{definition}

For $k=1$,  a $F_1$-overpartition is an overpartition in which the first occurrence of a number may be overlined. For $k=2$, there are nine $F_2$-overpartitions of $4$.
\[(4),(3,1),(\overline{3},{1}),(2,2),(\overline{2},{2}),(2,1,1),(2,\overline{1},1),(1,1,1,1),(\overline{1},1,1,1).
\]

Let $\mathcal{L}_k$ (resp. $\mathcal{F}_k$) be the set of all $L_k$-overpartitions (resp. $F_k$-overpartitions). We will show that $\mathcal{L}_k$ and $\mathcal{F}_k$ are separable overpartition classes and give the generating function for the overpartitions in $\mathcal{L}_k$ and $\mathcal{F}_k$ respectively.
\begin{theorem}\label{add-thm-eqn-L-F-k}
\begin{equation}\label{eqn-L-k}
\sum_{\pi\in\mathcal{L}_k}z^{\ell_o(\pi)}q^{|\pi|}
=1+\sum_{s=1}^k\sum_{m=1}^\infty\sum_{j=1}^m(z^{j-1}+z^j)\frac{q^{k{j\choose 2}+sj+k(m-j)}}{(q;q)_{k(m-1)+s}}{{m-1}\brack{j-1}}_k,
\end{equation}
and
\begin{equation}\label{eqn-F-k}
\begin{split}
\sum_{\pi\in\mathcal{F}_k}z^{\ell_o(\pi)}q^{|\pi|}
&=1+\sum_{s=1}^k\sum_{m=1}^\infty\sum_{j=1}^mz^{j-1}\frac{q^{k{j\choose 2}+sj+k(m-j)}}{(q;q)_{k(m-1)+s}}{{m-1}\brack{j-1}}_k\\
&\qquad+\sum_{m=1}^\infty\sum_{j=1}^mz^{j}\frac{q^{k{j\choose 2}+km}}{(q;q)_{km}}{{m-1}\brack{j-1}}_k,
\end{split}
\end{equation}
where ${A\brack B}_k$ is the $q$-binomial coefficient, or Gaussian polynomial for nonnegative integers $A$ and $B$ defined as follows:
\[{A\brack B}_k=\left\{\begin{array}{ll}\frac{(q^k;q^k)_A}{(q^k;q^k)_B(q^k;q^k)_{A-B}},&\text{if }A\geq B\geq 0,\\
0,&\text{otherwise.}
\end{array}
\right.\]
\end{theorem}

This article is organized as follows.  In section \ref{proof_2}, we first recall the definition of the conjugate of an overpartition and introduce the $(r+1)$-repeating size of an overpartition, and then give a proof of Theorem \ref{add-thm-gen-m-a-es}. We will show Theorems \ref{mex}, \ref{delta_1} and \ref{P_mex_1} in Sections 3. With the aid of separable overpartition classes, we will prove Theorem \ref{add-thm-eqn-L-F-k} in Section \ref{proof_3}.

\section{Proof of Theorem \ref{add-thm-gen-m-a-es}}\label{proof_2}

 The objective of this section is to give a proof of Theorem \ref{add-thm-gen-m-a-es}. Bear in mind that, in this section, an overpartition is defined as a partition in which the last occurrence of a number may be overlined. In Section \ref{proof_21}, we recall the definition of the conjugate of an overpartition and introduce the $(r+1)$-repeating size of an overpartition, and then we give the relation between the minimal (resp. maximal) excludant size and the largest (resp. smallest) $(r+1)$-repeating size.
We will show \eqref{gen-mes-r} and \eqref{gen-maes-r} in Section \ref{proof_22} and Section \ref{proof_23} respectively.

\subsection{Conjugate and $(r+1)$-repeating size}\label{proof_21}

 The conjugate of an overpartition was introduced by Corteel and Lovejoy \cite{Corteel-Lovejoy-2004}. Let $\pi=(\pi_1,\pi_2,\ldots,\pi_\ell)$ be an overpartition.  The conjugate $\pi'$ of $\pi$ is an overpartition such that
 \begin{itemize}
 \item[(1)] for $1\leq t<\ell$, $f_{\overline{t}}(\pi')+f_t(\pi')=|\pi_t|-|\pi_{t+1}|$;

 \item[(2)] $f_{\overline{\ell}}(\pi')+f_\ell(\pi')=|\pi_\ell|$;

 \item[(3)] for $1\leq t\leq \ell$,
 \[f_{\overline{t}}(\pi')=\left\{\begin{array}{ll}
1,&\text{if }\pi_t \text{ is overlined},\\
0,&\text{if }\pi_t \text{ is non-overlined};
\end{array}\right.\]

 \item[(4)] for $t>\ell$, $f_{\overline{t}}(\pi')=f_{t}(\pi')=0$.

 \end{itemize}

For example, the conjugate $\pi'$ of  overpartitions $\pi$ of $4$ are given in the following tables.

\begin{center}
\begin{tabular}{|c|c|c|c|c|c|c|c|c|}
  \hline
  $\pi$&$(4)$&$(\overline{4})$&$(3,1)$&$(\overline{3},1)$&$(3,\overline{1})$&$(\overline{3},\overline{1})$&$(2,2)$&$(2,\overline{2})$\\
  \hline
 $\pi'$&$(1,1,1,1)$&$(1,1,1,\overline{1})$&$(2,1,1)$&$({2},1,\overline{1})$&$(\overline{2},1,1)$&$(\overline{2},1,\overline{1})$&$(2,2)$&$(2,\overline{2})$\\
  \hline
\end{tabular}
\end{center}

\begin{center}
\begin{tabular}{|c|c|c|c|c|c|c|}
  \hline
  $\pi$&$(2,1,1)$&$(\overline{2},1,1)$&$(2,1,\overline{1})$&$(\overline{2},1,\overline{1})$&$(1,1,1,1)$&$(1,1,1,\overline{1})$\\
  \hline
  $\pi'$&$(3,1)$&$(3,\overline{1})$&$(\overline{3},1)$&$(\overline{3},\overline{1})$&$(4)$&$(\overline{4})$\\
  \hline
\end{tabular}
\end{center}

For $j\geq 1$, we say that $j$ is an $(r+1)$-repeating size of an overpartition $\pi$ if there are at least $r+1$ parts of size $j$ in $\pi$, that is, $j$ is an $(r+1)$-repeating size of  $\pi$ if and only if $f_{\overline{j}}(\pi)+f_j(\pi)\geq r+1$.
We assume that $0$ is an $(r+1)$-repeating size of $\pi$.

In terms of the conjugate of an overpartition, we can get the following two theorems, which give the relation between the minimal  excludant size and the largest $(r+1)$-repeating size and the relation between the  maximal excludant size and the smallest $(r+1)$-repeating size.

\begin{theorem}\label{relation-min}
For $j\geq 0$, $k\geq 1$ and $n\geq 0$, the number of overpartitions $\pi$ of $n$ such that ${mes}(\pi;r)=k$ and there are $j$ parts of size greater than $k$ in $\pi$ equals the number of overpartitions $\lambda$ of $n$ such that the largest $(r+1)$-repeating size of $\lambda$ is $j$ and there are $k-1$ parts of size greater than $j$ in $\lambda$.
\end{theorem}

For example, there are two overpartitions  of $4$ such that ${mes}(\pi;2)=1$ and there is one part of size greater than $1$, which are $(4)$ and $(\overline{4})$. There are two overpartitions  of $4$ such that the largest $3$-repeating size is $1$ and there are no parts of size greater than $1$, which are $(1,1,1,1)$ and $(1,1,1,\overline{1})$.

For another example, there are four overpartitions  of $4$ such that ${mes}(\pi;2)=4$ and there are no parts of size greater than $4$, which are $(3,1)$, $(\overline{3},1)$, $(3,\overline{1})$ and $(\overline{3},\overline{1})$. There are four overpartitions  of $4$ such that the largest $3$-repeating size is $0$ and there are three parts of size greater than $0$, which are $(2,1,1)$, $(2,1,\overline{1})$, $(\overline{2},1,1)$ and $(\overline{2},1,\overline{1})$.

\begin{theorem}\label{relation-max}
For $j\geq 1$, $k\geq 1$ and $n\geq 0$, the number of overpartitions $\pi$ of $n$ such that ${maes}(\pi;r)=k$ and there are $j$ parts of size greater than $k$ in $\pi$ equals the number of overpartitions $\lambda$ of $n$ such that the smallest positive $(r+1)$-repeating size of $\lambda$ is $j$ and there are $k+1$ parts of size greater than or equal to $j$ in $\lambda$.
\end{theorem}

For example, there are two overpartitions  of $6$ such that ${maes}(\pi;2)=2$ and there are two parts of size greater than $2$ in $\pi$, which are $(3,3)$ and $(3,\overline{3})$. There are two overpartitions  of $6$ such that the smallest positive $3$-repeating size is $2$ and there are three parts of size greater than or equal to $2$, which are $(2,2,2)$ and $(2,2,\overline{2})$.

We will give the proofs of \eqref{gen-mes-r} and \eqref{gen-maes-r} with the aid of Theorem \ref{relation-min} and Theorem \ref{relation-max} respectively.

\subsection{Proof of \eqref{gen-mes-r}}\label{proof_22}

In this subsection, we will give the proof of \eqref{gen-mes-r}. Before doing this, we first show that
 \begin{equation}\label{proof-mes-1}
 \sum_{j=0}^{\infty}\frac{q^{(r+1)j}(-1;q)_j}{(q;q)_j}\omega(j+1;\infty)=\frac{(-q;q)_\infty}{(q;q)_\infty},
 \end{equation}
 and for $n\geq 1$,
  \begin{equation}\label{proof-mes-2}
 \sum_{j=0}^{n-1}\frac{q^{(r+1)j}(-1;q)_j}{(q;q)_j}\omega(j+1;\infty)=\frac{(-q;q)_\infty}{(q;q)_\infty}-\omega(1;n-1)\frac{(-q^n;q)_\infty}{(q^n;q)_\infty}+\omega(1;\infty).
 \end{equation}

{\noindent \bf Proofs of \eqref{proof-mes-1} and \eqref{proof-mes-2}.} Obviously, the right hand side of \eqref{proof-mes-1} is the generating function for the overpartitions in $\overline{\mathcal{P}}$, that is,
     \begin{equation}\label{proof-mes-1-1}
     \sum_{\pi\in\overline{\mathcal{P}}}q^{|\pi|}=\frac{(-q;q)_\infty}{(q;q)_\infty}.
     \end{equation}

It is clear that for $j\geq 0$,
\[\frac{q^{(r+1)j}(-1;q)_j}{(q;q)_j}\omega(j+1;\infty)\]
is the generating function for the overpartitions with the largest $(r+1)$-repeating size being $j$. So, the left hand side of \eqref{proof-mes-1} is also the generating function for the overpartitions in $\overline{\mathcal{P}}$, that is,
     \begin{equation*}\label{proof-mes-1-2}
     \sum_{\pi\in\overline{\mathcal{P}}}q^{|\pi|}=\sum_{j=0}^{\infty}\frac{q^{(r+1)j}(-1;q)_j}{(q;q)_j}\omega(j+1;\infty).
     \end{equation*}
     Combining with \eqref{proof-mes-1-1}, arrive at \eqref{proof-mes-1}.

        For $n\geq 1$, the left hand side of \eqref{proof-mes-2} is the generating function for the overpartitions with the largest $(r+1)$-repeating size not exceeding $n-1$. It is easy to see that
        \[\omega(1;n-1)\frac{(-q^n;q)_\infty}{(q^n;q)_\infty}\]
        is the generating function for the overpartitions with the largest $(r+1)$-repeating size being $0$ or greater than $n-1$, and
        \[\omega(1;\infty)\]
        is the generating function for the overpartitions with the largest $(r+1)$-repeating size being $0$. Recalling that
        \[\frac{(-q;q)_\infty}{(q;q)_\infty}\]
        is the generating function for the overpartitions in $\overline{\mathcal{P}}$, we derive that  the right hand side of \eqref{proof-mes-2} is also the generating function for the overpartitions with the largest $(r+1)$-repeating size not exceeding $n-1$, and so \eqref{proof-mes-2} is valid. The proof is complete.  \qed

        Now, we are in a position to give the proof of  \eqref{gen-mes-r}.

        {\noindent \bf Proof of  \eqref{gen-mes-r}.} It follows from Theorem \ref{relation-min} that
     \begin{equation*}
     \sum_{\pi\in\overline{\mathcal{P}}}z^{mes(\pi;r)}q^{|\pi|}=\sum_{j=0}^\infty z\frac{q^{(r+1)j}(-1;q)_{j}}{(q;q)_j}\prod_{n=j+1}^{\infty}(1+2zq^n+2z^2q^{2n}+\cdots+2z^rq^{rn}).
     \end{equation*}
  Then, we get
\begin{align*}
&\quad\sum_{n=0}^\infty\sigma_r{mes}(n)q^n\\
&=\sum_{\pi\in\overline{\mathcal{P}}}mes(\pi;r)q^{|\pi|}\\
&=\frac{\partial}{\partial z}\left(\sum_{\pi\in\overline{\mathcal{P}}}z^{mes(\pi;r)}q^{|\pi|}\right)\bigg|_{z=1}\\
&=\sum_{j=0}^\infty \frac{q^{(r+1)j}(-1;q)_{j}}{(q;q)_j}\frac{\partial}{\partial z}\left(z\prod_{n=j+1}^{\infty}(1+2zq^n+2z^2q^{2n}+\cdots+2z^rq^{rn})\right)\bigg|_{z=1}\\
&=\sum_{j=0}^\infty \frac{q^{(r+1)j}(-1;q)_{j}}{(q;q)_j}\left(z\prod_{n=j+1}^{\infty}(1+2zq^n+2z^2q^{2n}+\cdots+2z^rq^{rn})\right)\bigg|_{z=1}\\
&\qquad\qquad\times\left(\frac{1}{z}+\sum_{n=j+1}^{\infty}\frac{2q^n+4zq^{2n}+\cdots+2rz^{r-1}q^{rn}}{1+2zq^n+2z^2q^{2n}+\cdots+2z^rq^{rn}}\right)\bigg|_{z=1}\\
&=\sum_{j=0}^\infty \frac{q^{(r+1)j}(-1;q)_{j}}{(q;q)_j}\omega(j+1;\infty)\left(1+\sum_{n=j+1}^{\infty}\frac{2(q^n+2q^{2n}+\cdots+rq^{rn})}{\omega(n)}\right)\\
&=\sum_{j=0}^\infty \frac{q^{(r+1)j}(-1;q)_{j}}{(q;q)_j}\omega(j+1;\infty)\\
&\quad+\sum_{j=0}^\infty\frac{q^{(r+1)j}(-1;q)_{j}}{(q;q)_j}\omega(j+1;\infty) \sum_{n=j+1}^{\infty}\frac{2(q^n+2q^{2n}+\cdots+rq^{rn})}{\omega(n)}\\
&=\sum_{j=0}^\infty \frac{q^{(r+1)j}(-1;q)_{j}}{(q;q)_j}\omega(j+1;\infty)\\
&\quad+\sum_{n=1}^{\infty}\frac{2(q^n+2q^{2n}+\cdots+rq^{rn})}{\omega(n)}
\sum_{j=0}^{n-1}\frac{q^{(r+1)j}(-1;q)_{j}}{(q;q)_j}\omega(j+1;\infty)\\
&=\frac{(-q;q)_{\infty}}{(q;q)_{\infty}}\\
&\qquad+\sum_{n=1}^{\infty}\frac{2(q^n+2q^{2n}+\cdots+rq^{rn})}{\omega(n)}\left(\frac{(-q;q)_{\infty}}{(q;q)_{\infty}}-\omega(1;n-1)\frac{(-q^n;q)_{\infty}}{(q^n;q)_{\infty}}+\omega(1;\infty)\right),
\end{align*}
where the final equation follows from \eqref{proof-mes-1} and \eqref{proof-mes-2}. The proof is complete. \qed

\subsection{Proof of \eqref{gen-maes-r}}\label{proof_23}

This subsection is devoted to presenting the proof of \eqref{gen-maes-r}. To do this, we are required to prove that
\begin{equation}\label{proof-maes-1}
 \sum_{j=1}^\infty 2q^{(r+1)j}\omega(1;j-1)\frac{(-q^{j+1};q)_\infty}{(q^j;q)_\infty}=\frac{(-q;q)_\infty}{(q;q)_\infty}-\omega(1;\infty),
 \end{equation}
and for $m\geq 1$,
\begin{equation}\label{proof-maes-2}
\sum_{j=1}^{m}2q^{(r+1)j}\omega(1;j-1)\frac{(-q^{j+1};q)_\infty}{(q^j;q)_\infty}=\frac{(-q;q)_\infty}{(q;q)_\infty}-\omega(1;m)\frac{(-q^{m+1};q)_\infty}{(q^{m+1};q)_\infty}.
 \end{equation}

 {\noindent \bf Proofs of \eqref{proof-maes-1} and \eqref{proof-maes-2}.} It is clear that for $j\geq 1$,
 \[2q^{(r+1)j}\omega(1;j-1)\frac{(-q^{j+1};q)_\infty}{(q^j;q)_\infty}\]
 is the generating function for the overpartitions with the smallest positive $(r+1)$-repeating size being $j$.
So, the left hand side of \eqref{proof-maes-1} is the generating function for the overpartitions having positive $(r+1)$-repeating size and the left hand side of \eqref{proof-maes-2} is the generating function for the overpartitions with smallest positive $(r+1)$-repeating size not exceeding $m$ for $m\geq 1$.

 Recalling  that
 \[\frac{(-q;q)_\infty}{(q;q)_\infty}\]
is the generating function for the overpartitions in $\overline{\mathcal{P}}$, and
 \[\omega(1;\infty)\]
is the generating function for the overpartitions with the largest $(r+1)$-repeating size being $0$, we deduce that the right hand side of \eqref{proof-maes-1} is also the generating function for the overpartitions having positive $(r+1)$-repeating size. So, \eqref{proof-maes-1} is valid.

It is clear that for $m\geq 1$, the right hand side of \eqref{proof-maes-2} is also the generating function for the overpartitions with smallest positive $(r+1)$-repeating size not exceeding $m$, since
\[\omega(1;m)\frac{(-q^{m+1};q)_\infty}{(q^{m+1};q)_\infty}\]
is the generating function for the overpartitions whose $(r+1)$-repeating sizes are $0$ or greater than $m$. So, \eqref{proof-maes-2} holds for $m\geq 1$. This completes the proof. \qed

We now turn to the proof of of \eqref{gen-maes-r}.

{\noindent \bf Proof of of \eqref{gen-maes-r}.} It follows from Theorem \ref{relation-max} that
     \begin{equation*}
     \sum_{\pi\in\overline{\mathcal{P}},\ maes(\pi;r)>0}z^{maes(\pi;r)}q^{|\pi|}=\sum_{j=1}^\infty 2z^rq^{(r+1)j}\omega(1;j-1)\frac{(-zq^{j+1};q)_\infty}{(zq^j;q)_\infty}.
     \end{equation*}
     Then, we get
\begin{align}
&\quad\sum_{n=0}^\infty\sigma_r{maes}(n)q^n\nonumber\\
&=\sum_{\pi\in\overline{\mathcal{P}}}maes(\pi;r)q^{|\pi|}\nonumber\\
&=\frac{\partial}{\partial z}\left(\sum_{\pi\in\overline{\mathcal{P}},\ maes(\pi;r)>0}z^{maes(\pi;r)}q^{|\pi|}\right)\bigg|_{z=1}\nonumber\\
&=\sum_{j=1}^\infty 2q^{(r+1)j}\omega(1;j-1)\frac{\partial}{\partial z}\left(z^r\frac{(-zq^{j+1};q)_\infty}{(zq^j;q)_\infty}\right)\bigg|_{z=1}\nonumber\\
&=\sum_{j=1}^\infty 2q^{(r+1)j}\omega(1;j-1)\left(z^r\frac{(-zq^{j+1};q)_\infty}{(zq^j;q)_\infty}\right)\bigg|_{z=1}\left(\frac{r}{z}+\sum_{n=j+1}^\infty\frac{q^n }{1+zq^n}-\sum_{n=j}^\infty\frac{-q^n}{1-zq^n}\right)\bigg|_{z=1}\nonumber\\
&=\sum_{j=1}^\infty 2q^{(r+1)j}\omega(1;j-1)\frac{(-q^{j+1};q)_\infty}{(q^j;q)_\infty}\left(r+\sum_{n=j+1}^\infty\frac{q^n }{1+q^n}+\sum_{n=j}^\infty\frac{q^n}{1-q^n}\right)\label{proof-maes-3}.
\end{align}

Using \eqref{proof-maes-2}, we have
\begin{align*}
&\quad\sum_{j=1}^\infty 2q^{(r+1)j}\omega(1;j-1)\frac{(-q^{j+1};q)_\infty}{(q^j;q)_\infty}\sum_{n=j+1}^\infty\frac{q^n}{1+q^n}\\
&=\sum_{n=2}^\infty\frac{q^n}{1+q^n}\sum_{j=1}^{n-1} 2q^{(r+1)j}\omega(1;j-1)\frac{(-q^{j+1};q)_\infty}{(q^j;q)_\infty}\\
&=\sum_{n=2}^\infty\frac{q^n}{1+q^n}\left(\frac{(-q;q)_\infty}{(q;q)_\infty}-\omega(1;n-1)\frac{(-q^{n};q)_\infty}{(q^{n};q)_\infty}\right)\\
&=\sum_{n=1}^\infty\frac{q^n}{1+q^n}\left(\frac{(-q;q)_\infty}{(q;q)_\infty}-\omega(1;n-1)\frac{(-q^{n};q)_\infty}{(q^{n};q)_\infty}\right)\\
&=\frac{(-q;q)_\infty}{(q;q)_\infty}\sum_{n=1}^\infty\frac{q^n}{1+q^n}-\sum_{n=1}^\infty q^n\omega(1;n-1)\frac{(-q^{n+1};q)_\infty}{(q^{n};q)_\infty},
\end{align*}
and
\begin{align*}
&\quad\sum_{j=1}^\infty 2q^{(r+1)j}\omega(1;j-1)\frac{(-q^{j+1};q)_\infty}{(q^j;q)_\infty}\sum_{n=j}^\infty\frac{q^n}{1-q^n}\\
&=\sum_{n=1}^\infty\frac{q^n}{1-q^n}\sum_{j=1}^{n} 2q^{(r+1)j}\omega(1;j-1)\frac{(-q^{j+1};q)_\infty}{(q^j;q)_\infty}\\
&=\sum_{n=1}^\infty\frac{q^n}{1-q^n}\left(\frac{(-q;q)_\infty}{(q;q)_\infty}-\omega(1;n)\frac{(-q^{n+1};q)_\infty}{(q^{n+1};q)_\infty}\right)\\
&=\frac{(-q;q)_\infty}{(q;q)_\infty}\sum_{n=1}^\infty\frac{q^n}{1-q^n}-\sum_{n=1}^\infty q^n\omega(1;n)\frac{(-q^{n+1};q)_\infty}{(q^{n};q)_\infty}.
\end{align*}
Combining with \eqref{proof-maes-1} and \eqref{proof-maes-3}, we get
\begin{align*}
&\quad\sum_{n=0}^\infty\sigma_r{maes}(n)q^n\\
&=r\sum_{j=1}^\infty 2q^{(r+1)j}\omega(1;j-1)\frac{(-q^{j+1};q)_\infty}{(q^j;q)_\infty}+\sum_{j=1}^\infty 2q^{(r+1)j}\omega(1;j-1)\frac{(-q^{j+1};q)_\infty}{(q^j;q)_\infty}\sum_{n=j+1}^\infty\frac{q^n}{1+q^n}\\
&\qquad+\sum_{j=1}^\infty 2q^{(r+1)j}\omega(1;j-1)\frac{(-q^{j+1};q)_\infty}{(q^j;q)_\infty}\sum_{n=j}^\infty\frac{q^n}{1-q^n}\\
&=r\left(\frac{(-q;q)_\infty}{(q;q)_\infty}-\omega(1;\infty)\right)+\frac{(-q;q)_\infty}{(q;q)_\infty}\sum_{n=1}^\infty\frac{q^n}{1+q^n}-\sum_{n=1}^\infty q^n\omega(1;n-1)\frac{(-q^{n+1};q)_\infty}{(q^{n};q)_\infty}\\
&\qquad+\frac{(-q;q)_\infty}{(q;q)_\infty}\sum_{n=1}^\infty\frac{q^n}{1-q^n}-\sum_{n=1}^\infty q^n\omega(1;n)\frac{(-q^{n+1};q)_\infty}{(q^{n};q)_\infty}\\
&=r\left(\frac{(-q;q)_\infty}{(q;q)_\infty}-\omega(1;\infty)\right)+\frac{(-q;q)_\infty}{(q;q)_\infty}\sum_{n=1}^\infty\left(\frac{q^n}{1+q^n}+\frac{q^n}{1-q^n}\right)\\
&\qquad-\sum_{n=1}^\infty q^n\left(\omega(1;n-1)+\omega(1;n)\right)\frac{(-q^{n+1};q)_\infty}{(q^{n};q)_\infty}\\
&=r\left(\frac{(-q;q)_{\infty}}{(q;q)_{\infty}}-\omega(1;\infty)\right)+\frac{(-q;q)_{\infty}}{(q;q)_{\infty}}\sum_{n=1}^\infty\frac{2q^n}{1-q^{2n}}\\
&\qquad\qquad-\sum_{n=1}^\infty q^n\omega(1;n-1)(1+w(n))\frac{(-q^{n+1};q)_{\infty}}{(q^n;q)_{\infty}}.
\end{align*}
We arrive at \eqref{gen-maes-r}, and thus the proof is complete.  \qed

\section{Proof of Theorems \ref{mex}, \ref{delta_1} and \ref{P_mex_1}}
In this section, an overpartition is defined as a partition in which the last occurrence of a number may be overlined. We will show  Theorems \ref{mex}, \ref{delta_1} and \ref{P_mex_1} in Sections \ref{submex}, \ref{subdelta_1} and \ref{subP_mex_1} respectively.

\subsection{Proof of Theorem \ref{mex}}\label{submex}
The objective of this subsection is to give a proof of Theorem \ref{mex}. To do this, we first show the following lemma, which will be used in the proof of Theorem \ref{delta_1}.
\begin{lemma}\label{add-lem-1}
For $s>m\geq 1$ and $n\geq 0$, let $\overline{A}_{m,s}(n)$ denote the number of overpartitions of $n$ with minimal excludant $m$ and second
minimal excludant $s$. Then, we have
\begin{equation*}\label{Ars_1}
\sum_{n=0}^\infty\overline{A}_{m,s}(n)q^n=\frac{(-q;q)_\infty}{(q;q)_\infty}q^{\binom{s}{2}-m}2^{s-2}\frac{(1-q^m)(1-q^s)}{(-q;q)_s}.
\end{equation*}
\end{lemma}

\pf Clearly, we have
\begin{align*}
\sum_{n=0}^\infty\overline{A}_{m,s}(n)q^n&=2^{m-1}q^{1+2+\cdots+(m-1)}2^{s-m-1}q^{(m+1)+(m+2)+\cdots+(s-1)}\frac{(1-q^m)(1-q^s)}{(q;q)_\infty}(-q^{s+1};q)_\infty\\
&=\frac{(-q;q)_\infty}{(q;q)_\infty}q^{\binom{s}{2}-m}2^{s-2}\frac{(1-q^m)(1-q^s)}{(-q;q)_s}.
\end{align*}
The proof is complete.   \qed

 Then, we proceed to show Theorem \ref{mex}.

 {\noindent\bf Proof of Theorem \ref{mex}.} By Lemma \ref{add-lem-1}, we have
\begin{equation}\label{Ars_2}
\begin{split}
&\quad\sum_{n=0}^\infty\sum_{s=2}^\infty\sum_{m=1}^{s-1}\overline{A}_{m,s}(n)\omega^sq^n\\
&=\sum_{s=2}^\infty\sum_{m=1}^{s-1}\omega^s \sum_{n=0}^\infty\overline{A}_{m,s}(n)q^n\\
&=\frac{(-q;q)_\infty}{(q;q)_\infty}\sum_{s=2}^\infty\sum_{m=1}^{s-1}\omega^sq^{\binom{s}{2}-m}2^{s-2}\frac{(1-q^m)(1-q^s)}{(-q;q)_s}\\
&=\frac{(-q;q)_\infty}{(q;q)_\infty}\sum_{s=2}^\infty\omega^s\frac{q^{\binom{s}{2}}2^{s-2}(1-q^s)}{(-q;q)_s}\sum_{m=1}^{s-1}\frac{1-q^m}{q^m}\\
&=\frac{(-q;q)_\infty}{(q;q)_\infty}\sum_{s=2}^\infty\omega^s\frac{q^{\binom{s}{2}}2^{s-2}(1-q^s)}{(-q;q)_s}\left(\frac{q^{-(s-1)}(1-q^{s-1})}{1-q}-(s-1)\right)\\
&=\frac{(-q;q)_\infty}{(1-q)(q;q)_\infty}\sum_{s=2}^\infty\omega^s\frac{q^{\binom{s-1}{2}}2^{s-2}(1-q^s)}{(-q;q)_s}\left(1-q^{s-1}-(s-1)(1-q)q^{s-1}\right).
\end{split}
\end{equation}

It is clear that for $n\geq 0$ and $s\geq 2$,
\[\sum_{m=1}^{s-1}\overline{A}_{m,s}(n)\]
is the number of overpartitions of $n$ with second minimal excludant $s$. This implies that
\[{\sigma}_2\overline{mex}(n)=\sum_{s=2}^{\infty}s\sum_{m=1}^{s-1}\overline{A}_{m,s}(n).\]
Combining with \eqref{Ars_2}, we get
\begin{align*}
\sum_{n=0}^\infty{\sigma}_2\overline{mex}(n)q^n&=\sum_{n=0}^\infty\sum_{s=2}^{\infty}s\sum_{m=1}^{s-1}\overline{A}_{m,s}(n)q^n\\
&=\frac{\partial}{\partial\omega}\left(\sum_{n=0}^\infty\sum_{s=2}^\infty\sum_{m=1}^{s-1}\overline{A}_{m,s}(n)\omega^sq^n\right)\bigg|_{\omega=1}\\
&=\frac{(-q;q)_\infty}{(1-q)(q;q)_\infty}\sum_{s=2}^\infty \frac{sq^{\binom{s-1}{2}}2^{s-2}(1-q^s)}{(-q;q)_s}\left(1-q^{s-1}-(s-1)(1-q)q^{s-1}\right).
\end{align*}

So, it remains to show
\begin{equation*}\label{add-new-eqn-001}
\sum_{s=2}^\infty\frac{sq^{\binom{s-1}{2}}2^{s-2}(1-q^s)}{(-q;q)_s}\left(1-q^{s-1}-(s-1)(1-q)q^{s-1}\right)=1+(1-q)\sum_{s=0}^\infty\frac{(2-s)q^{\binom{s+1}{2}}2^{s-1}}{(-q;q)_s}.
\end{equation*}
To do this, we just need to prove that
\begin{equation}\label{add-new-eqn-001}
\begin{split}
&\quad\sum_{s=2}^\infty \frac{sq^{\binom{s-1}{2}}2^{s-2}(1-q^s)}{(-q;q)_s}\left(1-q^{s-1}-(s-1)(1-q)q^{s-1}\right)\\
&=\sum_{s=1}^\infty \frac{q^{\binom{s}{2}}2^{s-1}}{(-q;q)_s}\left(1-q^s+(1-q)(s-2sq^{s})\right),
\end{split}
\end{equation}

\begin{equation}\label{add-new-eqn-002}
\sum_{s=1}^\infty \frac{q^{\binom{s}{2}}2^{s-1}}{(-q;q)_s}(1-q^s)=1,
\end{equation}
and
\begin{equation}\label{add-new-eqn-003}
\sum_{s=1}^\infty \frac{q^{\binom{s}{2}}2^{s-1}}{(-q;q)_s}(s-2sq^{s})=\sum_{s=0}^\infty\frac{(2-s)q^{\binom{s+1}{2}}2^{s-1}}{(-q;q)_s}.
\end{equation}

{\noindent \bf Proof of \eqref{add-new-eqn-001}.} Clearly, we have
\begin{align*}
&\quad\sum_{s=2}^\infty \frac{sq^{\binom{s-1}{2}}2^{s-2}(1-q^s)}{(-q;q)_s}\left(1-q^{s-1}-(s-1)(1-q)q^{s-1}\right)\\
&=\sum_{s=2}^\infty \frac{sq^{\binom{s-1}{2}}2^{s-2}(1+q^s-2q^s)}{(-q;q)_s}\left(1-q^{s-1}-(s-1)(1-q)q^{s-1}\right)\\
&=\sum_{s=1}^\infty \frac{(s+1)q^{\binom{s}{2}}2^{s-1}}{(-q;q)_{s}}\left(1-q^{s}-s(1-q)q^{s}\right)\\
&\qquad-\sum_{s=2}^\infty \frac{sq^{\binom{s-1}{2}+s}2^{s-1}}{(-q;q)_s}\left(1-q^{s-1}-(s-1)(1-q)q^{s-1}\right)\\
&=\frac{2(1-q)^2}{1+q}+\sum_{s=2}^\infty \frac{q^{\binom{s}{2}}2^{s-1}}{(-q;q)_s}\\
&\qquad\qquad\qquad\quad\times\left((s+1)(1-q^{s}-s(1-q)q^{s})-sq(1-q^{s-1}-(s-1)(1-q)q^{s-1})\right)\\
&=\sum_{s=1}^\infty \frac{q^{\binom{s}{2}}2^{s-1}}{(-q;q)_s}\left(1-q^s+(1-q)(s-2sq^{s})\right).
\end{align*}
We arrive at \eqref{add-new-eqn-001}.

{\noindent \bf Proof of \eqref{add-new-eqn-002}.} Obviously, we have
\begin{align*}
\sum_{s=1}^\infty \frac{q^{\binom{s}{2}}2^{s-1}}{(-q;q)_s}(1-q^s)
&=\sum_{s=1}^\infty \frac{q^{\binom{s}{2}}2^{s-1}}{(-q;q)_s}(1+q^s-2q^s)\\
&=\sum_{s=1}^\infty \frac{q^{\binom{s}{2}}2^{s-1}}{(-q;q)_{s-1}}-\sum_{s=1}^\infty \frac{q^{\binom{s+1}{2}}2^{s}}{(-q;q)_s}=1,
\end{align*}
and so \eqref{add-new-eqn-002} is valid.

{\noindent \bf Proof of \eqref{add-new-eqn-003}.} It is easy to see that
\begin{align*}
\sum_{s=1}^\infty \frac{sq^{\binom{s}{2}}2^{s-1}}{(-q;q)_s}&=\sum_{s=1}^\infty \frac{sq^{\binom{s}{2}}2^{s-1}}{(-q;q)_s}(1+q^s-q^s)\\
&=\sum_{s=1}^\infty \frac{sq^{\binom{s}{2}}2^{s-1}}{(-q;q)_{s-1}}-\sum_{s=1}^\infty \frac{sq^{\binom{s+1}{2}}2^{s-1}}{(-q;q)_s}\\
&=\sum_{s=0}^\infty \frac{(s+1)q^{\binom{s+1}{2}}2^{s}}{(-q;q)_s}-\sum_{s=0}^\infty \frac{sq^{\binom{s+1}{2}}2^{s-1}}{(-q;q)_s}.
\end{align*}
So, we have
\begin{align*}
\sum_{s=1}^\infty \frac{q^{\binom{s}{2}}2^{s-1}}{(-q;q)_s}(s-2sq^{s})&=\sum_{s=0}^\infty \frac{(s+1)q^{\binom{s+1}{2}}2^{s}}{(-q;q)_s}-\sum_{s=0}^\infty \frac{sq^{\binom{s+1}{2}}2^{s-1}}{(-q;q)_s}-\sum_{s=1}^\infty \frac{sq^{\binom{s+1}{2}}2^{s}}{(-q;q)_s}\\
&=\sum_{s=0}^\infty \frac{q^{\binom{s+1}{2}}2^{s}}{(-q;q)_s}-\sum_{s=0}^\infty \frac{sq^{\binom{s+1}{2}}2^{s-1}}{(-q;q)_s}.
\end{align*}
This implies that \eqref{add-new-eqn-003} holds, and thus the proof is complete. \qed

\subsection{Proof of Theorem \ref{delta_1}}\label{subdelta_1}

In this subsection, we aim to give a proof of Theorem \ref{delta_1} with the aid of Lemma \ref{add-lem-1}. For $n\geq 0$ and $t\geq 1$, it follows from the definition of $\overline{\Delta}_t(n)$ that
\begin{equation*}
\overline{\Delta}_t(n)=\sum_{m=1}^\infty\overline{A}_{m,m+t}(n).
\end{equation*}

It yields that
\begin{equation*}
\sum_{n=0}^{\infty}\overline{\Delta}_t(n)q^n=\sum_{n=0}^{\infty}\sum_{m=1}^\infty\overline{A}_{m,m+t}(n)q^n=\sum_{m=1}^\infty\left(\sum_{n=0}^{\infty}\overline{A}_{m,m+t}(n)q^n\right).
\end{equation*}
Combining with Lemma \ref{add-lem-1}, we can get
\begin{align*}
&\quad\sum_{n=0}^{\infty}\overline{\Delta}_t(n)q^n\\
&=\frac{(-q;q)_\infty}{(q;q)_\infty}\sum_{m=1}^\infty q^{\binom{m+t}{2}-m}2^{m+t-2}\frac{(1-q^m)(1-q^{m+t})}{(-q;q)_{m+t}}\\
&=\frac{(-q;q)_\infty}{(q;q)_\infty}\sum_{m=1}^\infty \frac{q^{\binom{m+t}{2}-m}2^{m+t-2}(1-q^m)}{(-q;q)_{m+t}}(1+q^{m+t}-2q^{m+t})\\
&=\frac{(-q;q)_\infty}{(q;q)_\infty}\left(\sum_{m=1}^\infty \frac{q^{\binom{m+t}{2}-m}2^{m+t-2}(1-q^m)}{(-q;q)_{m+t-1}}-
\sum_{m=1}^\infty \frac{q^{\binom{m+t+1}{2}-m}2^{m+t-1}(1-q^m)}{(-q;q)_{m+t}}\right)\\
&=\frac{(-q;q)_\infty}{(q;q)_\infty}\left(\sum_{m=0}^\infty \frac{q^{\binom{m+t+1}{2}-m-1}2^{m+t-1}(1-q^{m+1})}{(-q;q)_{m+t}}-
\sum_{m=0}^\infty \frac{q^{\binom{m+t+1}{2}-m}2^{m+t-1}(1-q^m)}{(-q;q)_{m+t}}\right)\\
&=\frac{(-q;q)_\infty}{(q;q)_\infty}\sum_{m=0}^\infty \frac{q^{\binom{m+t+1}{2}-m-1}2^{m+t-1}}{(-q;q)_{m+t}}\left(1-q^{m+1}-q(1-q^m)\right)\\
&=\frac{(-q;q)_\infty}{(q;q)_\infty}\sum_{m=0}^\infty \frac{q^{\binom{m+t+1}{2}-m-1}2^{m+t-1}}{(-q;q)_{m+t}}(1-q)\\
&=q^{t-1}(1-q)\frac{(-q;q)_\infty}{(q;q)_\infty}\sum_{m=0}^\infty \frac{q^{\binom{m+t+1}{2}-m-t}2^{m+t-1}}{(-q;q)_{m+t}}\\
&=(q^{t-1}-q^t)\frac{(-q;q)_\infty}{(q;q)_\infty}\sum_{m=0}^\infty\frac{q^{\binom{m+t}{2}}2^{m+t-1}}{(-q;q)_{m+t}}.
\end{align*}
We arrive at \eqref{delta-2}. This completes the proof of Theorem \ref{delta_1}.   \qed

\subsection{Proof of Theorem \ref{P_mex_1}}\label{subP_mex_1}

In this subsection, we will give a proof of Theorem \ref{P_mex_1}. For $k\geq 0$, let $\pi$ be an overpartition such that $\overline{mex}(\pi)=k+1$ and the length of the mex sequence of $\pi$ is at least $r$. Then,
\begin{itemize}
\item[(1)] there exist parts of size $1,2,\ldots,k$ in $\pi$;
\item[(2)] there are no parts of size $k+1,k+2,\ldots,k+r$ in $\pi$.
\end{itemize}
This implies that
\begin{align}
\sum_{n=0}^\infty\overline{p}_r^{mex}(n)q^n&=\sum_{k=0}^{\infty}\frac{2^{k}q^{1+2+\cdots+k}}{(q;q)_k}\frac{(-q^{k+r+1};q)_\infty}{(q^{k+r+1};q)_\infty}\nonumber\\
&=\sum_{k=0}^\infty \frac{q^{\binom{k+1}{2}}}{(q;q)_k}2^k\frac{(-q^{k+r+1};q)_\infty}{(q^{k+r+1};q)_\infty}.\label{proof11-1}
\end{align}

Setting $a=-1$ and $t=q^{k+r+1}$ in the $q$-binomial theorem \cite[Theorem 2.1]{Andrews-1976}:
\begin{equation*}
\sum_{m=0}^\infty \frac{(a;q)_{m}}{(q;q)_{m}}t^{m}=\frac{{(at;q)}_{\infty}}{(t;q)_{\infty}},
\end{equation*}
we can get
\begin{equation}\label{proof11-2}
\sum_{m=0}^\infty\frac{(-1;q)_m}{(q;q)_m}q^{(k+r+1)m}=\frac{(-q^{k+r+1};q)_\infty}{(q^{k+r+1};q)_\infty}.
\end{equation}

Substituting \eqref{proof11-2} into \eqref{proof11-1}, we get
\begin{align*}
\sum_{n=0}^\infty\overline{p}_r^{mex}(n)q^n&=\sum_{k=0}^\infty \frac{q^{\binom{k+1}{2}}}{(q;q)_k}2^k\sum_{m=0}^\infty\frac{(-1;q)_m}{(q;q)_m}q^{(k+r+1)m}\\
&=\sum_{m=0}^\infty\frac{(-1;q)_m}{(q;q)_m}q^{(r+1)m}\sum_{k=0}^\infty \frac{q^{\binom{k+1}{2}}}{(q;q)_k}(2q^{m})^k\\
&=\sum_{m=0}^\infty\frac{(-1;q)_m}{(q;q)_m}q^{(r+1)m}(-2q^{m+1};q)_\infty,
\end{align*}
where the final equation follows by letting $t=2q^{m+1}$ the following identity due to Euler \cite{Euler-1748} (see also \cite[(2.2.6)]{Andrews-1976}):
\begin{equation*}
\sum_{k=0}^\infty\frac{t^kq^{{k}\choose 2}}{(q;q)_k}=(-t;q)_\infty.
\end{equation*}
This completes the proof. \qed

\section{Proof of Theorem \ref{add-thm-eqn-L-F-k}}\label{proof_3}

The objective of this section is to give a proof of Theorem \ref{add-thm-eqn-L-F-k}. More precisely, we will show that $\mathcal{L}_k$ and $\mathcal{F}_k$ are separable overpartition classes and then give the proofs of  \eqref{eqn-L-k} and \eqref{eqn-F-k} in Section \ref{proof_31} and Section \ref{proof_32} respectively.

\subsection{Proof of \eqref{eqn-L-k}}\label{proof_31}

In this subsection, an overpartition is defined as a partition in which the last occurrence of a number may be overlined. The objective of this subsection is to give a proof \eqref{eqn-L-k}. We first show that $\mathcal{L}_k$ is a separable overpartition class. To do this, we are required to find the basis of $\mathcal{L}_k$, which involves the following set.

\begin{definition}
For $m\geq 1$, let $\mathcal{BL}_k(m)$ be the set of overpartitions $\lambda=(\lambda_1,\lambda_2,\ldots,\lambda_m)$ in $\mathcal{L}_k$ such that
\begin{itemize}
\item[(1)] $\lambda_m=\overline{1}$  or   $1$;

\item[(2)] for $1\leq i<m$, $|\lambda_i|\leq |\lambda_{i+1}|+1$ with strict inequality if $\lambda_i$ is non-overlined.
\end{itemize}
\end{definition}

For $m\geq 1$, let $\lambda$ be an overpartition in $\mathcal{BL}_k(m)$. For $1\leq i<m$, assume that  $|\lambda_{i+1}|=t$. If $m-i\equiv0\pmod k$, then we have $\lambda_i=t$ or $\overline{t+1}$. If $m-i\not\equiv0\pmod k$, then $\lambda_i$ is non-overlined, and so $\lambda_i=t$. Then, we obtain that for $m\geq 1$ and $1\leq s\leq k$, the number of overpartitions in $\mathcal{BL}_k(k(m-1)+s)$ is $2^m$. For example, the number of overpartitions in $\mathcal{BL}_2(5)$ is $2^3=8$.
\[(1,1,1,1,1),(\overline{2},1,1,1,1),(2,2,\overline{2},1,1),(\overline{3},2,\overline{2},1,1),\]
\[(1,1,1,1,\overline{1}),(\overline{2},1,1,1,\overline{1}),(2,2,\overline{2},1,\overline{1}),(\overline{3},2,\overline{2},1,\overline{1}).\]
\begin{lemma}
 $\mathcal{L}_k$ is a separable overpartition class.
\end{lemma}

\pf  Set \[\mathcal{BL}_k=\bigcup_{m=1}^\infty\mathcal{BL}_k(m).\]
 Obviously, $\mathcal{BL}_k$ is the basis of $\mathcal{L}_k$. This completes the proof.   \qed

For $m,j\geq 1$, let $\mathcal{BL}_k(m,j)$ (resp. $\mathcal{BL}_k(m,\overline{j}))$ be the set of overpartitions in $\mathcal{BL}_k(m)$ with the largest part being $j$ (resp. $\overline{j}$).
\begin{prop}\label{m>j}
Assume that $j\geq 1$,
\begin{itemize}
\item[(1)] for $2\leq s\leq k$ and $m\geq 1$, or $s=1$ and $m=1$, if $\mathcal{BL}_k(k(m-1)+s,j)$ is nonempty, then we have $m\geq j$;

\item[(2)] for $1\leq s\leq k$ and $m\geq 1$, if $\mathcal{BL}_k(k(m-1)+s,\overline{j})$ is nonempty, then we have $s = 1$ and $m\geq j$.

    \item[(3)] for $m\geq 2$, if $\mathcal{BL}_k(k(m-1)+1,j)$ is nonempty, then we have $m>j$.
    \end{itemize}
\end{prop}

\pf For $1\leq s\leq k$ and $m,j\geq 1$, assume that $\lambda$ is an overpartition in $\mathcal{BL}_k(k(m-1)+s)$ with the largest part of size $j$, then we see that
\[f_{\overline{t}}(\lambda)=1\text{ for }2\leq t\leq j,\]
\[f_{\overline{t}}(\lambda)+f_t(\lambda)\geq k\text{ and }f_{\overline{t}}(\lambda)+f_t(\lambda)\equiv 0\pmod{k}\text{ for }1\leq t<j,\]
and
\[f_{\overline{j}}(\lambda)+f_j(\lambda)\geq s\text{ and }f_{\overline{j}}(\lambda)+f_j(\lambda)\equiv s\pmod{k}.\]

So, we get
\begin{equation}\label{length-lambda-1}
\ell(\lambda)=\sum_{t=1}^{j}(f_{\overline{t}}(\lambda)+f_t(\lambda))\geq k(j-1)+s.
\end{equation}

Under the condition that $\lambda$ is an overpartition in $\mathcal{BL}_k(k(m-1)+s)$, we have $\ell(\lambda)=k(m-1)+s.$ Combining with \eqref{length-lambda-1}, we arrive at $m\geq j$, and so (1) is valid.

If the largest part of $\lambda$ is $\overline{j}$, then we have $\lambda_1=\overline{j}$. In this case, we have $s-1\equiv\ell(\lambda)-1\equiv0\pmod{k}$, and so $s=1$. Hence, (2) is verified.

Now, we proceed to show (3). Assume that $m\geq 2$ and $\lambda$ is an overpartition in $\mathcal{BL}_k(k(m-1)+1,j)$. Clearly, we have $\ell(\lambda)=k(m-1)+1$ and $f_{j}(\lambda)>0$. We consider the following two cases.

Case 1: $j=1$. In such case, it follows from the assumption $m\geq 2$ that $m>j$.

Case 2: $j\geq 2$. In such case, we have $f_{\overline{j}}(\lambda)=1$. Note that $f_{j}(\lambda)>0$, so we have $f_{\overline{j}}(\lambda)+f_j(\lambda)>1$. Then, we get
\[
k(m-1)+1=\ell(\lambda)=\sum_{t=1}^{j}(f_{\overline{t}}(\lambda)+f_t(\lambda))>k(j-1)+1,
\]
which yields $m>j$. This completes the proof of (3).    \qed

For $m,j\geq 1$, define
\[GL_k(m,j)= \sum_{\lambda\in\mathcal{BL}_k(m,j)}z^{\ell_o(\lambda)}q^{|\lambda|},\]
and
\[GL_k(m,\overline{j})=\sum_{\lambda\in\mathcal{BL}_k(m,\overline{j})}z^{\ell_o(\lambda)}q^{|\lambda|}.\]

We find that in order to prove \eqref{eqn-L-k}, it suffices to show the following theorem.
\begin{theorem}\label{equivalent-L-k}
For $1\leq s\leq k$ and $m\geq j\geq 1$,
\begin{equation}\label{eqn-eqv-L-K}
GL_k(k(m-1)+s,j)+GL_k(k(m-1)+s,\overline{j})=(z^{j-1}+z^j)q^{k{j\choose 2}+sj+k(m-j)}{{m-1}\brack{j-1}}_k.
\end{equation}
\end{theorem}

To prove Theorem \ref{equivalent-L-k}, we need the following lemmas.
\begin{lemma}\label{proof-equivalent-lemma-1}
For $2\leq s\leq k$  and $m\geq j\geq 2 $,
\[GL_k(k(m-1)+s,j)=(z^{j-1}+z^j)q^{k{j\choose 2}+sj+k(m-j)}{{m-1}\brack{j-1}}_k.\]
\end{lemma}
\pf Let $\lambda$ be an overpartition in $\mathcal{BL}_k(k(m-1)+s,j)$. If we remove the parts $\overline2,(k-1)$'s $2,\ldots,\overline{j-1},(k-1)$'s $j-1,\overline j,(s-1)$'s $j$ and $k$'s $1$ (resp. $\overline1,(k-1)$'s $1$) from $\lambda$ if $\overline{1}$ does not occur in $\lambda$ (resp. $\overline{1}$ occurs in $\lambda$), then we get a partition $\mu$ such that
\begin{itemize}
\item[(1)] $\ell(\mu)=k(m-j)$ and the parts in $\mu$ do not exceed $j$;
  \item[(2)]  $k|f_{t}(\mu)$ for $1\leq t\leq j$.

\end{itemize}
The process above to get $\mu$ could be run in reverse. Clearly, the generating function for the partitions $\mu$ satisfying the conditions (1) and (2) is
\[q^{k(m-j)}{{m-1}\brack{j-1}}_k.\]

So, we get
\[GL_k(k(m-1)+s,j)=(z^{j-1}+z^j)q^{k{j\choose 2}+sj}\cdot q^{k(m-j)}{{m-1}\brack{j-1}}_k.\]
This completes the proof.  \qed

\begin{lemma}\label{proof-equivalent-lemma-2}
For $m\geq j\geq 2$,
\[GL_k(k(m-1)+1,\overline{j})=(z^{j-1}+z^j)q^{k{j\choose 2}+j+k(m-j)}{{m-2}\brack{j-2}}_k.\]
\end{lemma}

\pf Let $\lambda$ be an overpartition in $\mathcal{BL}_k(k(m-1)+1,\overline{j})$. If we remove the parts $\overline2,(k-1)$'s $2,\ldots,\overline{j-1},(k-1)$'s $j-1,\overline j$ and $k$'s $1$ (resp. $\overline1,(k-1)$'s $1$) from $\lambda$ if $\overline{1}$ does not occur in $\lambda$ (resp. $\overline{1}$ occurs in $\lambda$), then we get a partition $\mu$ such that
\begin{itemize}
\item[(1)] $\ell(\mu)=k(m-j)$ and the parts in $\mu$ do not exceed $j-1$;
  \item[(2)]  $k|f_{t}(\mu)$ for $1\leq t\leq j-1$.

\end{itemize}
The process above to get $\mu$ could be run in reverse. Clearly, the generating function for the partitions $\mu$ satisfying the conditions (1) and (2) is
\[q^{k(m-j)}{{m-2}\brack{j-2}}_k.\]

So, we get
\[GL_k(k(m-1)+1,\overline{j})=(z^{j-1}+z^j)q^{k{j\choose 2}+j}\cdot q^{k(m-j)}{{m-2}\brack{j-2}}_k,\]
and thus the proof is complete.  \qed

\begin{lemma}\label{proof-equivalent-lemma-3}
For $m>j\geq 2$,
 \[GL_k(k(m-1)+1,j)=(z^{j-1}+z^j)q^{k{j\choose 2}+j+k(m-1)}{{m-2}\brack{j-1}}_k.\]
\end{lemma}
\pf Let $\lambda$ be an overpartition in $\mathcal{BL}_k(k(m-1)+1,j)$. If we remove the parts $\overline2,(k-1)$'s $2,\ldots,\overline{j-1},(k-1)$'s $j-1,\overline j,k$'s $j$ and $k$'s $1$ (resp. $\overline1,(k-1)$'s $1$) from $\lambda$ if $\overline{1}$ does not occur in $\lambda$ (resp. $\overline{1}$ occurs in $\lambda$), then we get a partition $\mu$ such that
\begin{itemize}
\item[(1)] $\ell(\mu)=k(m-j-1)$ and the parts in $\mu$ do not exceed $j$;
  \item[(2)]  $k|f_{t}(\mu)$ for $1\leq t\leq j$.

\end{itemize}
The process above to get $\mu$ could be run in reverse. Clearly, the generating function for the partitions $\mu$ satisfying the conditions (1) and (2) is
\[q^{k(m-j-1)}{{m-2}\brack{j-1}}_k.\]

So, we get
\[GL_k(k(m-1)+1,j)=(z^{j-1}+z^j)q^{k{j\choose 2}+(k+1)j}\cdot q^{k(m-j-1)}{{m-2}\brack{j-1}}_k.\]
We complete the proof.  \qed

Now, we are in a position to give a proof of Theorem \ref{equivalent-L-k}.

{\noindent \bf Proof of Theorem \ref{equivalent-L-k}.} For $j=1$, it is clear that
\[GL_k(m,\overline{1})=\left\{\begin{array}{ll}zq, & m=1, \\
 0, & m\geq 2,
 \end{array}\right.\]
and
\[GL_k(m,1)=\left\{\begin{array}{ll}q, & m=1, \\
 (1+z)q^m, & m\geq 2.
 \end{array}\right.\]
 This implies that for $1\leq s\leq k$ and $m\geq 1$,
\[
GL_k(k(m-1)+s,1)+GL_k(k(m-1)+s,\overline{1})=(1+z)q^{k(m-1)+s},\]
which agrees with \eqref{eqn-eqv-L-K} for $j=1$.

For $2\leq s\leq k$  and $m\geq j\geq 2 $, note that $GL_k(k(m-1)+s,\overline{j})=0$, by Lemma \ref{proof-equivalent-lemma-1}, we have
\begin{align*}
&\quad GL_k(k(m-1)+s,j)+GL_k(k(m-1)+s,\overline{j})\\
&=GL_k(k(m-1)+s,j)\\
&=(z^{j-1}+z^j)q^{k{j\choose 2}+sj+k(m-j)}{{m-1}\brack{j-1}}_k.
\end{align*}

For $m=j\geq 2$, it follows from the condition (3) in Proposition \ref{m>j} that $GL_k(k(m-1)+1,{j})=0$. Using Lemma \ref{proof-equivalent-lemma-2}, we get
\begin{align*}
&\quad GL_k(k(m-1)+1,j)+GL_k(k(m-1)+1,\overline{j})\\
&=GL_k(k(m-1)+1,\overline{j})\\
&=(z^{j-1}+z^j)q^{k{j\choose 2}+j+k(m-j)},
\end{align*}
which is \eqref{eqn-eqv-L-K} for $s=1$ and $m=j\geq 2$.

For $s=1$ and $m>j\geq 2$, appealing to Lemmas  \ref{proof-equivalent-lemma-2} and \ref{proof-equivalent-lemma-3}, we get
\begin{align*}
&\quad GL_k(k(m-1)+1,j)+GL_k(k(m-1)+1,\overline{j})\\
&=(z^{j-1}+z^j)q^{k{j\choose 2}+j+k(m-j)}{{m-2}\brack{j-2}}_k+(z^{j-1}+z^j)q^{k{j\choose 2}+j+k(m-1)}{{m-2}\brack{j-1}}_k\\
&=(z^{j-1}+z^j)q^{k{j\choose 2}+j+k(m-j)}\left( {{m-2}\brack{j-2}}_k+q^{k(j-1)}{{m-2}\brack{j-1}}_k\right)\\
&=(z^{j-1}+z^j)q^{k{j\choose 2}+j+k(m-j)}{{m-1}\brack{j-1}}_k,
\end{align*}
where the final equation follows from the standard recurrence for the $q$-binomial coefficients \cite[(3.3.4)]{Andrews-1976}:
\begin{equation*}\label{bin-new-r-1}
{A\brack B}_k={{A-1}\brack{B-1}}_k+q^{kB}{{A-1}\brack{B}}_k.
\end{equation*}
This completes the proof.  \qed

We conclude this subsection with the following property.
\begin{prop}\label{prop-bl-D}
For $1\leq s\leq k$ and $n,j\geq 1$, let ${BL}'_{k,s}(n,j)$ (resp. ${BL}''_{k,s}(n,j-1)$) be the number of overpartitions $\lambda$ of $n$ in $\mathcal{BL}_k$ such that $\ell(\lambda)\equiv s\pmod{k}$, $\ell_o(\lambda)=j$ (resp. $\ell_o(\lambda)=j-1$) and the smallest part of $\lambda$ is $\overline{1}$ (resp. $1$), and let ${D}_{k,s}(n,j)$ be the number of partitions of $n$ with $j$ distinct parts which are congruent to $s$ modulo $k$. Then, we have
\[{BL}'_{k,s}(n,j)={BL}''_{k,s}(n,j-1)={D}_{k,s}(n,j).\]
\end{prop}

We will give an analytic proof and a combinatorial proof of Proposition  \ref{prop-bl-D}.

{\noindent \bf Analytic proof of Proposition  \ref{prop-bl-D}.}  Clearly, we have
\begin{align*}
&\quad\sum_{n=1}^\infty\sum_{j=1}^\infty {BL}'_{k,s}(n,j)z^{j}q^{n}+\sum_{n=1}^\infty\sum_{j=1}^\infty {BL}''_{k,s}(n,j-1)z^{j-1}q^{n}\\
&=\sum_{m=1}^\infty \sum_{j=1}^m(GL_k(k(m-1)+s,j)+GL_k(k(m-1)+s,\overline{j})).
\end{align*}

In light of Theorem \ref{equivalent-L-k}, we get
\begin{align}
&\quad \sum_{n=1}^\infty\sum_{j=1}^\infty {BL}'_{k,s}(n,j)z^{j}q^{n}+\sum_{n=1}^\infty\sum_{j=1}^\infty {BL}''_{k,s}(n,j-1)z^{j-1}q^{n}\nonumber\\
&=\sum_{m=1}^\infty \sum_{j=1}^m(z^{j-1}+z^j)q^{k{j\choose 2}+sj+k(m-j)}{{m-1}\brack{j-1}}_k\nonumber\\
&=\sum_{j=1}^\infty (z^{j-1}+z^j)q^{k{j\choose 2}+sj}\sum_{m=j}^\infty q^{k(m-j)}{{m-1}\brack{j-1}}_k.\label{bin-q-q-k}
\end{align}

For $m\geq j\geq 1$, we know that
\[q^{m-j}{{m-1}\brack{j-1}}_1\]
is the generating function for the partitions $\mu$ such that $\ell(\mu)\leq j$ and the largest part of $\mu$ is $m-j$. So, we obtain that
\[\sum_{m=j}^\infty q^{m-j}{{m-1}\brack{j-1}}_1\]
is the generating function for the partitions $\mu$ such that $\ell(\mu)\leq j$, that is,
\begin{equation}\label{eqn-brack-1}
\sum_{m=j}^\infty q^{m-j}{{m-1}\brack{j-1}}_1=\frac{1}{(q;q)_j}.
\end{equation}

Letting $q\rightarrow q^k$ in \eqref{eqn-brack-1}, we get
 \begin{equation}\label{eqn-brack-2}
\sum_{m=j}^\infty q^{k(m-j)}{{m-1}\brack{j-1}}_k=\frac{1}{(q^k;q^k)_j}.
\end{equation}

Substituting \eqref{eqn-brack-2} into \eqref{bin-q-q-k}, we have
\begin{align*}
&\quad \sum_{n=1}^\infty\sum_{j=1}^\infty {BL}'_{k,s}(n,j)z^{j}q^{n}+\sum_{n=1}^\infty\sum_{j=1}^\infty {BL}''_{k,s}(n,j-1)z^{j-1}q^{n}\\
&=\sum_{j=1}^\infty \frac{(z^{j-1}+z^j)q^{k{j\choose 2}+sj}}{(q^k;q^k)_j}\\
&=(1+z^{-1})\sum_{j=1}^\infty \frac{(zq^s)^jq^{k{j\choose 2}}}{(q^k;q^k)_j}\\
&=(1+z^{-1})((-zq^s;q^k)_\infty-1),
\end{align*}
where the final equation follows by letting $q\rightarrow q^k$ and $z\rightarrow zq^s$ in \cite[(2.2.6)]{Andrews-1976}:
\begin{equation}\label{Euler-2}
\sum_{j=0}^\infty\frac{z^jq^{{j}\choose 2}}{(q;q)_j}=(-z;q)_\infty.
\end{equation}

Note that
\begin{equation}\label{gen-d}
\sum_{n=1}^\infty\sum_{j=1}^\infty {D}_{k,s}(n,j)z^{j}q^{n}=(-zq^s;q^k)_\infty-1,
\end{equation}
so we have
\begin{equation}\label{proof-oooo-1}
\sum_{n=1}^\infty\sum_{j=1}^\infty {BL}'_{k,s}(n,j)z^{j}q^{n}+\sum_{n=1}^\infty\sum_{j=1}^\infty {BL}''_{k,s}(n,j-1)z^{j-1}q^{n}=(1+z^{-1})\sum_{n=1}^\infty\sum_{j=1}^\infty {D}_{k,s}(n,j)z^{j}q^{n}.
\end{equation}

For $n,j\geq 1$, let $\lambda$ be an overpartition counted by ${BL}'_{k,s}(n,j)$. If we replace the smallest part $\overline{1}$ of $\lambda$ by a non-overlined part $1$, then we get an overpartition enumerated by ${BL}''_{k,s}(n,j-1)$, and vice versa. This implies that
\[\sum_{n=1}^\infty\sum_{j=1}^\infty {BL}''_{k,s}(n,j-1)z^{j-1}q^{n}=z^{-1}\sum_{n=1}^\infty\sum_{j=1}^\infty {BL}'_{k,s}(n,j)z^{j}q^{n}.\]
Combining with \eqref{proof-oooo-1}, we get
\[\sum_{n=1}^\infty\sum_{j=1}^\infty {BL}'_{k,s}(n,j)z^{j}q^{n}=\sum_{n=1}^\infty\sum_{j=1}^\infty {D}_{k,s}(n,j)z^{j}q^{n},\]
and
\[\sum_{n=1}^\infty\sum_{j=1}^\infty {BL}''_{k,s}(n,j-1)z^{j-1}q^{n}=\sum_{n=1}^\infty\sum_{j=1}^\infty {D}_{k,s}(n,j)z^{j-1}q^{n}.\]
This completes the proof.  \qed

{\noindent \bf Combinatorial proof of Proposition  \ref{prop-bl-D}.} For $1\leq s\leq k$ and $n,j\geq 1$, let $\mathcal{BL}'_{k,s}(n,j)$ and $\mathcal{BL}''_{k,s}(n,j-1)$ be the set of overpartitions counted by ${BL}'_{k,s}(n,j)$ and ${BL}''_{k,s}(n,j-1)$ respectively. In the analytic proof of Proposition  \ref{prop-bl-D}, we have proved that there is a bijection between $\mathcal{BL}'_{k,s}(n,j)$ and $\mathcal{BL}''_{k,s}(n,j-1)$.
Let $\mathcal{D}_{k,s}(n,j)$ be the set of partitions enumerated by ${D}_{k,s}(n,j)$, we just need to establish a bijection between $\mathcal{BL}'_{k,s}(n,j)$ and $\mathcal{D}_{k,s}(n,j)$.

Let $\lambda$ be an overpartition in $\mathcal{BL}'_{k,s}(n,j)$. Assume that $|\lambda_1|=j$,  we first remove the parts $\overline 1,(k-1)$'s $1,\overline 2,(k-1)$'s $2,\ldots,\overline{j-1},(k-1)$'s $j-1,\overline j$ and $(s-1)$'s $j$  from $\lambda$, then we get a partition $\mu$ of $n-k{j\choose 2}-sj$ such that
\begin{itemize}
\item[(1)] the parts in $\mu$ do not exceed $j$;

  \item[(2)]  $k|f_{t}(\mu)$ for $1\leq t\leq j$.

\end{itemize}

Let $\mu'$ be the conjugate of $\mu$, which is a partition such that the $i$-th part $\mu'_i$ of $\mu'$ is the number of parts of $\mu$ greater than or equal to $i$.
 Then, $\mu'$ is a partition such that $\ell(\mu')\leq j$ and the parts in $\mu'$ are divisible by $k$. For $\ell(\mu')<i\leq j$, we assume that $\mu'_i=0$. Set
\[\nu=(\mu'_1+k(j-1)+s,\ldots,\mu'_{j-1}+k+s,\mu'_j+s).\]
Clearly, $\nu$ is a partition  in $\mathcal{D}_{k,s}(n,j)$.  Obviously, the
process above is reversible. The proof is complete.   \qed

\subsection{Proof of \eqref{eqn-F-k}}\label{proof_32}

In this subsection, an overpartition is defined as a partition in which the first occurrence of a number may be overlined. The objective of this subsection is to give a proof \eqref{eqn-F-k}. We first show that $\mathcal{F}_k$ is a separable overpartition class. To do this, we are required to find the basis of $\mathcal{F}_k$, which involves the following set.

\begin{definition}
For $m\geq 1$, let $\mathcal{BF}_k(m)$ be the set of overpartitions $\lambda=(\lambda_1,\lambda_2,\ldots,\lambda_m)$ in $\mathcal{F}_k$ such that

\begin{itemize}
\item[(1)] if $k=1$, then $\lambda_m=\overline{1}$ or $1$; if $k\geq2$, then $\lambda_m=1$

\item[(2)] for $1\leq i<m$, $|\lambda_i|\leq |\lambda_{i+1}|+1$ with strict inequality if $\lambda_{i+1}$ is non-overlined.
\end{itemize}
\end{definition}

For $m\geq 1$, let $\lambda$ be an overpartition in $\mathcal{BF}_k(m)$.
For $1\leq i<m$, we consider the following two cases.

Case 1: if $\lambda_{i+1}=t$, then we have $|\lambda_i|\geq t$ and $|\lambda_i|<t+1$, which implies that $|\lambda_i|=t$. Moreover, we have $\lambda_i=t$ or $\overline{t}$ if $m-i\equiv-1\pmod k$, and $\lambda_i=t$ otherwise.

Case 2: if $\lambda_{i+1}=\overline{t}$, then we have $|\lambda_i|>{t}$ and $|\lambda_i|\leq t+1$, which implies that $|\lambda_i|=t+1$. Moreover, we have $\lambda_i=t+1$ or $\overline{t+1}$ if $m-i\equiv-1\pmod k$, and $\lambda_i=t+1$ otherwise.

Then, we obtain that for $m\geq 1$, the number of overpartitions in $\mathcal{BF}_k(k(m-1)+s)$ is $2^{m-1}$ for $1\leq s<k$ and the number of overpartitions in $\mathcal{BF}_k(k(m-1)+s)$ is $2^{m}$ for $s=k$.  For example, the number of overpartitions in $\mathcal{BF}_2(5)$ is $2^{3-1}=4$.
\[(1,1,1,1,1),(2,\overline{1},1,1,1),(2,2,2,\overline{1},1),(3,\overline{2},2,\overline{1},1).\]

\begin{lemma}
 $\mathcal{F}_k$ is a separable overpartition class.
\end{lemma}

\pf  Set \[\mathcal{BF}_k=\bigcup_{m=1}^\infty\mathcal{BF}_k(m).\]
 Obviously, $\mathcal{BF}_k$ is the basis of $\mathcal{F}_k$. This completes the proof.   \qed

For $m,j\geq 1$, let $\mathcal{BF}_k(m,j)$ (resp. $\mathcal{BF}_k(m,\overline{j}))$ be the set of overpartitions in $\mathcal{BF}_k(m)$ with the largest part being $j$ (resp. $\overline{j}$).

\begin{prop}\label{m>j-f}
For $1\leq s\leq k$ and $m,j\geq1$,
\begin{itemize}
\item[(1)] if $\mathcal{BF}_k(k(m-1)+s,j)$ is nonempty, then we have $m\geq j$;

\item[(2)] if $\mathcal{BF}_k(k(m-1)+s,\overline{j})$ is nonempty, then we have $s=k$ and $m\geq j$.

\end{itemize}

\end{prop}

\pf By definition, we obtain that if $\mathcal{BF}_k(k(m-1)+s,\overline{j})$ is nonempty  then  $s=k$. Assume that $\lambda$ is an overpartition in $\mathcal{BF}_k(k(m-1)+s)$ with the largest part of size $j$, then we see that
\[f_{\overline{t}}(\lambda)+f_t(\lambda)\geq k\text{ and }f_{\overline{t}}(\lambda)+f_t(\lambda)\equiv 0\pmod{k}\text{ for }1\leq t<j,\]
and
\[f_{\overline{j}}(\lambda)+f_j(\lambda)\geq s\text{ and }f_{\overline{j}}(\lambda)+f_j(\lambda)\equiv s\pmod{k}.\]

So, we get
\begin{equation}\label{length-lambda-1-f}
\ell(\lambda)=\sum_{t=1}^{j}(f_{\overline{t}}(\lambda)+f_t(\lambda))\geq k(j-1)+s.
\end{equation}

Under the assumption that $\lambda$ is an overpartition in $\mathcal{BF}_k(k(m-1)+s)$, we have $\ell(\lambda)=k(m-1)+s.$ Combining with \eqref{length-lambda-1-f}, we arrive at $m\geq j$. The proof is complete.    \qed

For $m,j\geq 1$, define
\[GF_k(m,j)= \sum_{\lambda\in\mathcal{BF}_k(m,j)}z^{\ell_o(\lambda)}q^{|\lambda|},\]
and
\[GF_k(m,\overline{j})=\sum_{\lambda\in\mathcal{BF}_k(m,\overline{j})}z^{\ell_o(\lambda)}q^{|\lambda|}.\]

We find that in order to prove \eqref{eqn-F-k}, it suffices to show the following theorem.
\begin{theorem}\label{equivalent-F-k}
For $1\leq s\leq k$ and $m\geq j\geq 1$,
\begin{equation}\label{eqn-eqv-F-K-1}
GF_k(k(m-1)+s,j)=z^{j-1}q^{k{j\choose 2}+sj+k(m-j)}{{m-1}\brack{j-1}}_k,
\end{equation}
and for $m\geq j\geq 1$,
\begin{equation}\label{eqn-eqv-F-K-2}
GF_k(km,\overline{j})=z^jq^{k{j\choose 2}+km}{{m-1}\brack{j-1}}_k.
\end{equation}
\end{theorem}

\pf For $m\geq j\geq 1$, let $\lambda$ be an overpartition in $\mathcal{BF}_k(km,\overline{j})$.
If we replace the largest part $\overline{j}$ of $\lambda$ by a non-overlined part $j$, then we get an overpartition in $\mathcal{BF}_k(km,{j})$, and vice versa. This implies that
\[GF_k(km,\overline{j})=z GF_k(km,{j}).\]

We just need to show \eqref{eqn-eqv-F-K-1}. For $1\leq s\leq k$ and $m\geq j\geq 1$, let $\lambda$ be an overpartition in $\mathcal{BF}_k(k(m-1)+s,{j})$. If we remove the parts $(k-1)$'s $1,\overline{1},(k-1)$'s $2,\overline{2},\ldots,(k-1)$'s $j-1,\overline{j-1}$ and $s'$s $j$ from $\lambda$, then we get a partition $\mu$ such that
\begin{itemize}
\item[(1)] $\ell(\mu)=k(m-j)$ and the parts in $\mu$ do not exceed $j$;
  \item[(2)]  $k|f_{t}(\mu)$ for $1\leq t\leq j$.

\end{itemize}
The process above to get $\mu$ could be run in reverse. Clearly, the generating function for the partitions $\mu$ satisfying the conditions (1) and (2) is
\[q^{k(m-j)}{{m-1}\brack{j-1}}_k.\]

So, we get
\[GF_k(k(m-1)+s,j)=z^{j-1}q^{k{j\choose 2}+sj}\cdot q^{k(m-j)}{{m-1}\brack{j-1}}_k.\]
This completes the proof.  \qed

Finally, we give the following property.
\begin{prop}\label{prop-bF-D}
For $1\leq s\leq k$ and $n,j\geq 1$, let ${BF}'_{k,k}(n,j)$  be the number of overpartitions $\lambda$ of $n$ in $\mathcal{BF}_k$ such that $\ell(\lambda)\equiv 0\pmod{k}$, $\ell_o(\lambda)=j$   and the largest part of $\lambda$ is $\overline{j}$, let   ${BF}''_{k,s}(n,j-1)$ be the number of overpartitions $\lambda$ of $n$ in $\mathcal{BF}_k$ such that $\ell(\lambda)\equiv s\pmod{k}$, $\ell_o(\lambda)=j-1$ and the largest part of $\lambda$ is $j$, and let ${D}_{k,s}(n,j)$ be the number of partitions of $n$ with $j$ distinct parts which are congruent to $s$ modulo $k$. Then, we have
\[{BF}'_{k,k}(n,j)={D}_{k,k}(n,j)\text{ and }{BF}''_{k,s}(n,j-1)={D}_{k,s}(n,j).\]
\end{prop}

We will give an analytic proof and a combinatorial proof of Proposition  \ref{prop-bF-D}.

{\noindent \bf Analytic proof of Proposition  \ref{prop-bF-D}.}
In view of \eqref{eqn-brack-2}, \eqref{Euler-2},  \eqref{gen-d}, \eqref{eqn-eqv-F-K-1} and \eqref{eqn-eqv-F-K-2}, we get that
\begin{align*}
&\quad\sum_{n=1}^\infty\sum_{j=1}^\infty {BF}'_{k,k}(n,j)z^{j}q^{n}\\
&=\sum_{m=1}^\infty \sum_{j=1}^m GL_k(km,\overline{j})\\
&=\sum_{m=1}^\infty \sum_{j=1}^mz^jq^{k{j\choose 2}+km}{{m-1}\brack{j-1}}_k\\
&=\sum_{j=1}^\infty z^jq^{k{j\choose 2}+kj}\sum_{m=j}^\infty q^{k(m-j)}{{m-1}\brack{j-1}}_k\\
&=(-zq^k;q^k)_\infty-1\\
&=\sum_{n=1}^\infty\sum_{j=1}^\infty {D}_{k,k}(n,j)z^{j}q^{n},
\end{align*}
and for $1\leq s\leq k$,
\begin{align*}
&\quad\sum_{n=1}^\infty\sum_{j=1}^\infty {BF}''_{k,s}(n,j-1)z^{j-1}q^{n}\\
&=\sum_{m=1}^\infty \sum_{j=1}^m GF_k(k(m-1)+s,j)\\
&=\sum_{m=1}^\infty \sum_{j=1}^m z^{j-1}q^{k{j\choose 2}+sj+k(m-j)}{{m-1}\brack{j-1}}_k\\
&=z^{-1}\sum_{j=1}^\infty z^jq^{k{j\choose 2}+sj}\sum_{m=j}^\infty q^{k(m-j)}{{m-1}\brack{j-1}}_k\\
&=z^{-1}((-zq^s;q^k)_\infty-1)\\
&=\sum_{n=1}^\infty\sum_{j=1}^\infty {D}_{k,s}(n,j)z^{j-1}q^{n}.
\end{align*}
This completes the proof. \qed

{\noindent \bf Combinatorial proof of Proposition  \ref{prop-bF-D}.} For $1\leq s\leq k$ and $n,j\geq 1$, let $\mathcal{D}_{k,s}(n,j)$ be the set of partitions enumerated by ${D}_{k,s}(n,j)$, and
let $\mathcal{BF}'_{k,k}(n,j)$ and $\mathcal{BF}''_{k,s}(n,j-1)$ be the set of overpartitions counted by ${BF}'_{k,k}(n,j)$ and ${BF}''_{k,s}(n,j-1)$ respectively.
 In the  proof of Theorem \ref{equivalent-F-k}, we have proved that there is a bijection between $\mathcal{BF}'_{k,k}(n,j)$ and $\mathcal{BF}''_{k,k}(n,j-1)$. We just need to build a bijection between $\mathcal{BF}''_{k,s}(n,j-1)$ and $\mathcal{D}_{k,s}(n,j)$.

Let $\lambda$ be an overpartition in $\mathcal{BF}''_{k,s}(n,j-1)$. Assume that $\lambda_1=j$,  we remove the parts $(k-1)$'s $1,\overline{1},(k-1)$'s $2,\overline{2},\ldots,(k-1)$'s $j-1,\overline{j-1}$ and $s'$s $j$  from $\lambda$ and denote the resulting partition by $\mu$. With a similar argument in the combinatorial proof of Proposition  \ref{prop-bl-D}, we can get a partition
\[\nu=(\mu'_1+k(j-1)+s,\ldots,\mu'_{j-1}+k+s,\mu'_j+s)\in\mathcal{D}_{k,s}(n,j),\]
where   $\mu'$ is the conjugate of $\mu$.  Obviously, the
process above is reversible. The proof is complete.   \qed

\noindent{\bf Acknowledgments.} This work
was supported by   the National Natural Science Foundation of China (No. 12571356) and Sichuan Science and Technology Program (No. 2024NSFSC1394).

\end{document}